\newcommand{\Label}[1]{\label{#1}\hspace{.3cm}\fbox{\rm #1}\hspace{.3cm}}
\renewcommand{\Label}{\label} 
\newcommand{\abel}[1]{\marginpar{\ref{#1}=#1}}
\renewcommand{\abel}[1]{}  
\documentclass[11pt]{article}
 \usepackage{amsmath,amsfonts,amssymb}
 \usepackage{amsfonts}
 \usepackage{amssymb}
 \usepackage[dvips]{epsfig}
\begin{document}
\newcommand{\ls}[1]
   {\dimen0=\fontdimen6\the\font \lineskip=#1\dimen0
 \advance\lineskip.5\fontdimen5\the\font \advance\lineskip-\dimen0
 \lineskiplimit=.9\lineskip \baselineskip=\lineskip
 \advance\baselineskip\dimen0 \normallineskip\lineskip
 \normallineskiplimit\lineskiplimit \normalbaselineskip\baselineskip
 \ignorespaces }
 \newtheorem{theorem}{Theorem}[section]
 \newtheorem{conjecture}[theorem]{Conjecture}
 \newtheorem{corollary}[theorem]{Corollary}
 \newtheorem{lemma}[theorem]{Lemma}
 \newtheorem{assumption}[theorem]{Assumption}
 \newtheorem{proposition}[theorem]{Proposition}
 \newtheorem{definition}[theorem]{Definition}
 \newcommand{\remark}{\noindent {\bf Remark:\ }}
 \newcommand{\proof}{\noindent {\bf Proof:\ }}
 \newcommand{\pf}{\noindent \mbox{{\bf Proof}: }}
 \def\squarebox#1{\hbox to #1{\hfill\vbox to #1{\vfill}}}
 \newcommand{\qed}{\hspace*{\fill}
    \vbox{\hrule\hbox{\vrule\squarebox{.667em}\vrule}\hrule}\smallskip}
 \newcommand{\req}[1]{(\ref{#1})}
\newcommand{\Leb}{{\rm Leb\,}}
\newcommand{\al}{\alpha}
 \newcommand{\liminfn}{\liminf_{n \rightarrow \infty}}%
 \newcommand{\limsupk}{\limsup_{k \rightarrow \infty}}%
 \newcommand{\limsupn}{\limsup_{n \rightarrow \infty}}%
 \newcommand{\limsupe}{\limsup_{\eps \rightarrow \infty}}%
 \newcommand{\lip}{\langle}
 \newcommand{\rip}{\rangle}
 \newcommand{\lf}{\lfloor}
 \newcommand{\lc}{\lceil}
 \newcommand{\rc}{\rceil}
 \newcommand{\rf}{\rfloor}
 \newcommand{\uu}{\underline}
 \newcommand{\oo}{\overline}
 \newcommand{\won}{{\boldsymbol 1}}
 \newcommand{\chara}[2]{\mbox{${\bf 1}_{#1}(#2)$}}
 \newcommand{\supp}{{\operatorname{supp\ }}}
 \newcommand{\sign}{{\operatorname{sign\ }}}
 \newcommand{\one}{{\frac{1}{n}}}%
 \newcommand{\half}{{\frac{1}{2}}}%
 \newcommand{\ffrac}[2]{\left( \frac{#1}{#2} \right)}
 \newcommand{\bfrac}[2]{\frac{\mbox{\large{$#1$}}}{\mbox{\large{$#2$}}}}
 \newcommand{\dfr}{\displaystyle\frac}
 \newcommand{\hilbert}{\bigcirc\kern -0.8em{\rm\scriptstyle {H}\;}}
 \newcommand{\ga}{\gamma}
 \newcommand{\La}{\Lambda}
 \newcommand{\la}{\lambda}
 \newcommand{\lo}{\lambda^{(1)}}
 \newcommand{\eps}{\varepsilon}
 \newcommand{\om}{\omega}
 \newcommand{\Om}{\Omega}
 \newcommand{\Bla}{B_\la}
 \newcommand{\Ceps}{C^{(\eps)}}
 \newcommand{\Cepsm}{C^{(-\eps)}}
 \newcommand{\DCFn}{DCF${}^{(n)}$}
 \newcommand{\dnc}{\Delta_C^{(n)}(k)}
 \newcommand{\gn}{g_0^{(n)}}
 \newcommand{\Gzn}{\calG_z^{(n)}}
 \newcommand{\Kn}{K^{(n)}}
 \newcommand{\Knk}{K^{(n)^k}}
 \newcommand{\MN}{Murnaghan--Nakayama}
 \newcommand{\nga}{\left\lfloor n^\gamma\right\rfloor}
 \newcommand{\pins}{\pi^{(n)}_S}
 \newcommand{\Pn}{\calP_n}
 \newcommand{\Pnb}{\calP_{n,\beta}}
 \newcommand{\thelan}{\theta_\la^{(n)}}
 \newcommand{\thelank}{\theta_\la^{(n)^k}}
 \newcommand{\wmu}{\widehat{\mu}}
 \newcommand{\Xla}{\chi_{_\la}}
 \newcommand{\Xlast}{\chi_{_{\la^*}}}
 \newcommand{\bfcdot}{{\boldsymbol \cdot}}
 \newcommand{\calA}{{\mathcal A}}
 \newcommand{\calB}{{\mathcal B}}
 \newcommand{\calD}{{\mathcal D}}%
 \newcommand{\EE}{{\mathbb E}}%
 \newcommand{\calF}{{\mathcal F}}%
 \newcommand{\calG}{{\mathcal G}}%
 \newcommand{\calL}{{\mathcal L}}%
 \newcommand{\calM}{{\mathcal M}}%
 \newcommand{\NN}{{\mathbb N}}%
 \newcommand{\PP}{{\mathbb P}}%
 \newcommand{\calP}{{\mathcal P}}%
 \newcommand{\calT}{{\mathcal T}}%
 \newcommand{\ZZ}{{\mathbb Z}}   %
 \newcommand{\Nat}{\mbox{${\rm I\!N}$}}
 \newcommand{\Reals}{{\reals}}%
 \newcommand{\reals}{{\mathbb R}}%
 \newcommand{\Ups}{\Upsilon}
 \newcommand{\Ul}{\Ups_{_{\!\la}}}
 \newcommand{\lij}{\la_*^{(i,j)}}
 \newcommand{\aaa}{\mbox{\boldmath $a$}}
 \newcommand{\bbb}{\mbox{\boldmath $b$}}
 \newcommand{\saa}{\mbox{\boldmath \scriptsize{$a$}}}
 \newcommand{\sbb}{\mbox{\boldmath \scriptsize{$b$}}}
 \newcommand{\calC}{\mathcal C}
 \newcommand{\dab}{\delta_{\saa,\sbb}}
 \newcommand{\Cab}{C_{\saa,\sbb}}
 \newcommand{\Gabn}{G^{(n)}_{\saa,\sbb}}
\begin{centering}
{\large \bf
The Poisson-Dirichlet law is the unique
invariant distribution for uniform
split-merge transformations}\\[2em]
{\sc Persi Diaconis},\footnote{Department
of  Mathematics and Department of
Statistics,
Stanford University, Stanford, CA 94305, USA.}\
{\sc Eddy Mayer-Wolf},\footnote{Department of Mathematics,
Technion, Haifa 32000, Israel (email: emw@tx.technion.ac.il).
Partially supported by the
S. Faust research fund.}\\
{\sc Ofer Zeitouni},\footnote{Departments of Math.\ and of EE,
Technion, Haifa 32000, Israel, and Dept.\ of Math., University of Minnesota,
MN 55455 (email:
zeitouni@math.umn.edu).
Partially supported
by the fund for promotion of research at
the Technion,
and by a US-Israel
BSF grant.}\
and \
{\sc
Martin P.W.\ Zerner} \footnote{Department of Mathematics, Stanford
University, Stanford, CA 94305, U.S.A.
(email:\\ zerner@math.stanford.edu). Partially supported
by an Aly Kaufman Fellowship at the Technion.}\\[2em]
July 2, 2002. Revised January 21, 2003. Note added August 15, 2003.\\[2em] {\bf Dedicated to
the memory of Bob Brooks (1952--2002)}

\end{centering}

\numberwithin{equation}{section}
\begin{abstract}
 We consider a Markov chain on the space of (countable) partitions
 of the interval $[0,1]$, obtained first by size biased sampling
 twice (allowing repetitions)  and then merging the parts (if the
 sampled parts are distinct) or splitting the part uniformly (if
 the same part was sampled twice). We prove a conjecture of Vershik
 stating that the Poisson-Dirichlet law with parameter\ \,$\theta\!=\!1$
 \ is the unique invariant distribution for this Markov chain.
 Our proof uses a combination of probabilistic, combinatoric, and
 representation-theoretic arguments.

\end{abstract}

\vfill

\newpage

\ls{1}
 \section{Introduction}  \Label{intro}
Let $\Omega_1$ denote the space of (ordered) partitions of $[0,1]$, that
is
\[\Omega_1:=\left\{p\in [0,1]^{\Nat}: \ p_1\geq p_2\geq
... \geq 0,\ |p|_1=1\right\},\]
where $|x|_1=\sum_i|x_i|$ for any finite or countable sequence $(x_i)$.
By {\it size-biased} sampling according to a point $p\in \Omega_1$
we mean picking the $j$-th  part $p_j$ with probability $p_j$.
Our interest in this paper
is in the following Markov chain on $\Omega_1$,
which we call a
{\it continuous coagulation-fragmentation} process (CCF):
size-bias sample (with replacement) two parts from $p$. If the same
part was picked twice, split it (uniformly), and reorder the partition.
If different parts were picked, merge them, and reorder the partition.

 \noindent
 We denote by \DCFn \ ({\it discrete coagulation-fragmentation}) the Markov
 chain describing the evolution of the cycle lengths of permutations of
 $\{1,\ldots,n\}$ under random transpositions. The CCF process appears in a
 variety of contexts, but of particular relevance to us is its occurrence
 as a natural limit of \DCFn, when $n$ increases, see \cite{Tsilevich} for a
 discussion of this and its link with the space of ``virtual permutations''.

 \noindent
 For any $n\in\Nat$ denote
\[  \Pn :=\{\ell=(\ell_i)_{i\ge 1}
  \in\{0,1,\ldots,n\}^{\Nat}\ :\ \ell_1\geq \ell_2\geq\cdots\geq 0, \
|\ell|_1=n
\}\ \subset n\Omega_1\,.%
\]
(Elements in $\Pn$ may be thought of as being of length $n$; the
 remaining entries are necessarily zero).

\noindent
 A sequence $\ell\in\Pn$ is uniquely determined by its {\it type}\
 $(N_\ell(k)\!=\!\sharp\{i:\ell_i\!\!=\!\!k\})_{_{k=1}}^{^n}$, with
 $N_\ell\!=\!\sum_{_{k=1}}^{^n}\!\!N_\ell(k)$\ denoting  $\ell$'s total
 number of parts.

 \noindent
 The long-time behaviour of the \DCFn, viewed as an evolution in $\Pn$,
 is well understood. In particular, see e.g. \cite{diaconis}, it possesses a
 unique stationary distribution given by the Ewens formula:
 \begin{equation}
 \label{ewens}
 \pins(\ell)=\left(\prod_{k=1}^n k^{N_\ell(k)}N_\ell(k)!\right)^{-1}=
          \left(\prod_{i=1}^n \ell_i\prod_{k=1}^n N_\ell(k)!\right)^{-1}\,,
                       \hspace{1.5cm}\ell\in\Pn\,.
 \end{equation}

 \noindent
 It is well known, at least since \cite{Ki1,Kinew,VS},  that the measures
 $ \pins(n\cdot)$ on $\Omega_1$ converge weakly to the Poisson-Dirichlet
 distribution $\wmu_1$ with parameter $\theta\!=\!1$ (a precise definition
 of $\wmu_1$ is given below in Section \ref{cpreliminaries}).
 It has been shown in more than one way (cf. \cite{gnedin,pitman,Tsilevich})
 that $\wmu_1$ is invariant for the CCF transition.
 This fact, and hints coming from the theory of virtual permutations,
 led Vershik (see \cite{Tsilevich}) to
 \begin{conjecture}[Vershik]
 \label{conj-vershik} $\wmu_1$ is
 the  unique invariant distribution for the CCF.
 \end{conjecture}

 \noindent
 Our goal in this article is to prove Vershik's conjecture. A naive
 approach toward the proof would be to use the link with the \DCFn\
 and the fact that the latter converges to the distribution
 $\pins$ exponentially fast. However, the rate of that convergence
 deteriorates with $n$. To overcome this difficulty, our strategy
 consists of the following steps:
\begin{enumerate}
\item We provide a-priori estimates
(Proposition \ref{m1}) showing that every invariant
distribution for the CCF  leads to a good control on the
number of ``small parts''.
\item We couple the \DCFn \ and the CCF in such a way that whenever
they start from initial distributions  with such control on the
tails, the decoupling time is roughly $\sqrt{n}$ (Theorem
\ref{m2})
\item
For initial conditions as above, and for an appropriate class of
test functions, we show by using some harmonic analysis on the
symmetric group that the \DCFn\ achieves near equilibrium before
the decoupling time (Theorem \ref{the-xxx}).
\end{enumerate}
These steps are then combined in Theorem \ref{euni} to yield the
proof of Vershik's conjecture.

 \noindent
 Our work began from discussions with Bob Brooks on various models for
 ``random Riemann surfaces". Brooks and Makover \cite{brooksmakover1},
 \cite{brooksmakover0} studied Riemann surfaces via a dense set of
 ``Belyi surfaces" associated to three-regular graphs on $n$ vertices
 with an orientation at each vertex. Their construction gives a complete
 Riemann surface with finite area $\pi n$ for each graph. Uniformly
 choosing a random three-regular graph gives a probability distribution
 on Riemann surfaces, see \cite{gamburdmakover} for an accessible account
 of this model. Practical choice of a random three-regular graph is not
 so easy when $n$ is large. Brooks proposed a Markov chain method which
 involved splitting and joining cycles; investigating properties of his
 algorithm gave rise to the present paper.

 \noindent
 We next review some of the literature on this question. Tsilevich,
 in \cite{Tsilevich}, proves that $\wmu_1$ is the only
 CCF-invariant measure that is also invariant under additional
 symmetry conditions. Pitman, in \cite{pitman}, proves that
 $\wmu_1$ is the only CCF-invariant measure which is also invariant
 under size-biased sampling.  Related results appear in \cite{gnedin1}. 
In
 another direction,  it is shown in  \cite{MZZ} that $\wmu_1$ is
 the only CCF-invariant measure that is analytic in the sense that
 for any $k$, the law of an independently size-biased sample (with
 replacement) possesses an analytic density. Finally, Tsilevich in
 \cite{Tsilevich1} shows that the law of the CCF, initialized  at
 $p=(1,0,\ldots)$, and stopped  at a Binomial($n,1/2$) random time,
 converges to $\wmu_1$.

 \noindent
 We conclude this introduction by noting that in \cite{MZZ}, we
 have introduced a slightly more general model of split-merge transformations,
 by allowing either the split or the merge operations to be rejected
 with a certain probability. An invariant measure for these
 generalizations is the Poisson-Dirichlet law of parameter
 $\theta>0$. The discrete counterpart of this chain has been analyzed in
\cite[Section 4]{diaconishanlon}.  While it is plausible
 that the techniques of the
 current paper can be adapted to that setup using the results of \cite{diaconishanlon},
we do not pursue this
 generalization here.

\section{Continuous and Discrete Coagulation-Fragmentation}
  \Label{CCF}
\subsection{Preliminaries and CCF}  \Label{cpreliminaries}
 Given a topological space $W$, its Borel $\sigma$-algebra will be denoted by
 $\calB_W$, and the space of probability measures on $(W,\calB_W)$ by
 $\calM_1(W)$. By a slight abuse of notations, $\calM_1(V)$ will also be
 $\calM_1(W)$'s subspace of probability measures whose support is contained in a
 given closed subset $V$ of $W$. The total variation of a measure $\nu$ is
 denoted by $\|\nu\|_{\mbox{\scriptsize var}}$.
 \smallskip

 \noindent
We equip $\Omega_1$
with its relative $|\cdot|_{_1}$--topology which,
on $\Omega_1$, coincides with the
 weak (coordinatewise convergence) topology.

 \noindent
 On $\Omega_1$ we consider the Markov chain CCF in which two segments
 $p_i$ and $p_j$ of a given partition $p$ are size-bias
 sampled with replacement
 and then, if $i\!\ne\!j$ they merge into one of
length $p_i+p_j$ (coagulation),
 while if $i\!=\!j,\  p_i$ splits into two new parts $up_i,(1-u)p_i$ with
 $u\sim U[0,1]$ independent of all the rest (fragmentation). In either case the
 new partition is then rearranged nonincreasingly.
  \smallskip

 \noindent
 Recall
that the Poisson-Dirichlet
  law $\wmu_1$ is invariant for the CCF transition.
  Indeed, $\wmu_1$
itself has been defined in a variety of manners (\cite{arratia,Ki1})
  which are well known to be equivalent. Perhaps the simplest
is the GEM description
in which segments are successively and
uniformly removed from whatever remains of
 $[0,1]$, and then rearranged nonincreasingly.
Namely, let $Y_1\!\!=\!\!1$ and for
$n\!\in\!\Nat$ define\ $X_n\!=\!U_n\,Y_n,\ $
(the removed part at stage $n$) and
  $\ Y_{n+1}\!=\!Y_n\!-\!X_n$
(the remaining segment from which the $(n+1)$-th part
  is to be removed),\ where the $U_n$'s are
independent $U[0,1]$ variables. Since
 $Y_{n+1}\!=\!(1\!-\!U_n)Y_n$ it follows that $1\!\!-\!\!Y_{n+1}\!
  =\!\sum_{i=1}^n X_i$ increases
almost surely to $1$ as $n\to\infty$.
  The distribution on $\Omega_1$ of the nonincreasing
rearrangement $(p_i)_i$ of
 $(X_n)_n$ is called the Poisson--Dirichlet
law (with parameter $\theta=1$) and
 denoted $\wmu_1$.

 As has been mentioned in the Introduction,
it is the ultimate goal of this work to
 show that the Poisson-Dirichlet law is the
{\it only} CCF-invariant probability
  distribution. It will be crucial
for the main argument to establish in advance that
 any such invariant distribution does not put too
much weight on very small parts:

\begin{proposition}\Label{m1}
Let $\mu\in\calM_1(\Omega_1)$ be CCF-invariant.
Then
\begin{equation}\Label{hi}
\int\sum_{i\geq 1}p_i^\alpha\ d\mu<\infty\qquad\mbox{for all $\alpha>2/5.$}
\end{equation}
\end{proposition}\abel{m1}\abel{hi}
 The proof is deferred to the Appendix.

\subsection{DCF}   \Label{DCF}
 In this section we formally introduce the coagulation--fragmentation chain on
 the discrete version of $\Omega_1$, in which the partition points lie on a
 finite equidistant grid in $[0,1]$, or its equivalent state space $\Pn$, the
 set of integer partitions of a fixed $n\in\Nat$ defined in the Introduction.
 It will be helpful to view $\Pn$ as the conjugacy classes of the permutation
 group $S_n$.

 \noindent
 The \DCFn\ Markov chain on $\Pn$ is defined similarly to the CCF chain on
 $\Omega_1$. Identify each $\ell\in\Pn$ with a partition $\bigcup_i\!A_i$
 of $\{1,2,\ldots,n\}$,
where for each $i$, $\ell_i$ denotes the cardinality of $A_i$,
 and sample two independent integers $x,y$ uniformly from $\{1,\ldots,n\}$
 and without replacement, say $x\!\!\in\!\!A_i$ and $y\!\!\in\!\!A_j$.
 If $i\ne j$ replace $A_i$ and $A_j$ by $A_i\cup A_j$ while if $i=j$
 (in which case $\ell_i\ge 2$ since $x\ne y\in A_i$) replace $A_i$ by
 two of its subsets, consisting respectively of $A_i$'s $k$ smallest
 elements and of the $\ell_i\!-\!k$ remaining ones, where $k$ is
 uniformly sampled from $\{1,\ldots,\ell_i\!-\!1\}$ independently of $x$ and $y$.
 In either case relabel and rearrange the new $A_i$'s if necessary.

 \noindent
 The transition matrix $\Kn$ of \DCFn\ is described as follows:
To  split into or merge two parts of different sizes $j$ and $k$
($1\!\le\!j\!<\!k\!\le\!n$),
let $\ell,\ell'\!\in\!\Pn$ be such that
 $N_{\ell'}(j)\!=\!N_\ell(j)\!-\!1,\ \,N_{\ell'}(k)\!
 =\!N_\ell(k)\!-\!1,\ N_{\ell'}(j\!+\!k)\!
 =\!N_\ell(j\!+\!k)\!+\!1$\ \,and\ \,
 $N_{\ell'}(q)\!=\!N_\ell(q)$ for all $q\!\not\in\!\{j,k,j+k\}$. Then
 \begin{eqnarray} \Label{Kn}
   \Kn(\ell,\ell')&=&\frac{2jk}{n(n-1)}\,N_\ell(j)N_\ell(k)
                                \hspace{1.45cm}\mbox{merge}   \nonumber \\
   \Kn(\ell',\ell)&=&\frac{2(j+k)}{n(n-1)}\,N_{\ell'}(j+k)\hspace{1.7cm}
\mbox{split.}
 \end{eqnarray}
To  split into or merge two parts of the same size $j$ with $2\leq 2j\leq n$
let $\ell,\ell'\in \Pn$  and
$0\leq N_{\ell'}(j)=N_\ell(j)-2$, $N_{\ell'}(2j)=N_\ell(2j)+1$,
and
 $N_{\ell'}(q)\!=\!N_\ell(q)$ for all $q\!\not\in\!\{j,2j\}$. Then
\begin{eqnarray} \Label{Kn1}
   \Kn(\ell,\ell')&=&\frac{j^2}{n(n-1)}\,N_\ell(j)(N_\ell(j)-1)
                                \hspace{1.45cm}\mbox{merge}   \nonumber \\
   \Kn(\ell',\ell)&=&\frac{2j}{n(n-1)}\,N_{\ell'}(2j)\hspace{3.2cm}\mbox{split}
 \end{eqnarray}

\noindent
All other entries of the transition kernel are zero.
 \smallskip

 \noindent
 It is customary to think of the representation
 $\{1,2,\ldots,n\}=\bigcup\!A_i$ above as the notation for the conjugacy class
 of a permutation $\sigma\in S_n$. Seen this way, the \DCFn\ transition is
 nothing but the action of a random transposition on $S_n$'s conjugacy classes.
 Since the random transposition's unique stationary probability measure is the
 uniform law on $S_n$ (being a finite group convolution), one concludes that the
 \DCFn's unique stationary probability measure is the one induced on $S_n$'s
 conjugacy classes by the uniform law, namely (\ref{ewens}) (the Ewens sampling
 formula).
 In fact, \DCFn\ is reversible with respect to $\pins$, which can also be
 checked directly by using~(\ref{Kn}), (\ref{Kn1}) and~(\ref{ewens}) to verify the detailed
 balance equation $\Kn(\ell,\ell')\pins(\ell)=\Kn(\ell',\ell)\pins(\ell')$.
 \smallskip

 \section{Coupling of CCF and DCF}   \Label{coupling}
 In order to successfully approximate a CCF chain by \DCFn\ chains as $n\to\infty$ it
 is necessary to couple them on a common probability space.

\begin{theorem}   \Label{m2}
For all $\mu\in\calM_1(\Omega_1)$ and $\al<1/2$ satisfying
\begin{equation}\Label{m7}
\int\sum_{i\geq 1}p_i^\alpha\ d\mu<\infty,
\end{equation}
it is possible to define for all $n\geq 1$ a CCF Markov
chain $p(k)\ (k\geq 0)$ with initial distribution $\mu$ and a
$DCF^{(n)}$ Markov chain $\ell(k)\ (k\geq 0)$
 on the same probability space
with probability measure $Q_\mu^{(n)}$ and expectation $E_\mu^{(n)}$
in such a way that
\begin{eqnarray}
\lim_{n\to\infty}Q_\mu^{(n)}\left[N_{\ell(0)}\geq n^\beta\right]&=&0
\quad \mbox{for all $\al<\beta$ and}\Label{m8} \\ \Label{m9}
\lim_{n\to\infty} E_\mu^{(n)}\left[\left|p(\lfloor
n^\beta\rfloor)- \frac{\ell(\lfloor n^\beta\rfloor)}{n}\right|_1
\right]&=&0 \quad \mbox{for all $\beta<1/2$.}
\end{eqnarray}\abel{m8} \abel{m9}
\end{theorem}\abel{m2}\abel{m7}

\begin{proof} Fix $n\geq 1$.
We shall construct a Markov chain $(c_k,d_k,e_k)\ (k\geq 0)$ on
the state space
\begin{eqnarray*}
\Omega_{cde}^{(n)}&:=&\big\{(c,d,e)\ |\ c:[0,n)\to \ZZ\
\mbox{measurable},\quad \ZZ\backslash c[[0,n)]\ \mbox{infinite,}\\
&&\hspace*{18mm} d:\{1,\ldots,n\}\to \ZZ,\quad e\in\{0,1\}\big\}.
\end{eqnarray*}
Here $c$ and $d$ describe a continuous partition of $[0,n)$ and a
discrete partition of $\{1,\ldots,n\}$, respectively. The
interpretation of $c$ and $d$ in terms of elements of $\Omega_1$
and $\Pn$ is given by the  functions $\pi_c:\
\Omega_{cde}^{(n)}\to \Omega_{1}$ and $\pi_d:\
\Omega_{cde}^{(n)}\to \calP_n$, respectively, defined by
\begin{eqnarray}
\pi_c(c,d,e)\Label{france} &:=&{\rm sort\,
}\left(\left(\frac{\Leb(c^{-1}(\{i\}))}{n}\right)_{i\in \ZZ}
\right)\qquad\mbox{and}\\ \pi_d(c,d,e) &:=& {\rm sort\, }((\sharp
d^{-1}(\{i\}))_{i\in \ZZ}). \Label{senegal}
\end{eqnarray} \abel{france}\abel{senegal}
Here Leb denotes the Lebesgue measure and
${\rm sort\, }((x_i)_i)$ is the sequence obtained by
arranging the $x_i$'s in decreasing order, ignoring the 0's
  if there are infinitely many positive $x_i$'s. Thus two points $x, y\in
[0,n)$  belong to the same set in the partition of $[0,n)$ which
is described by $c$   iff $c(x)=c(y)$. Analogously,
$x,y\in\{1,\ldots,n\}$ belong to the same set in the partition of
$\{1,\ldots,n\}$ described by $d$  iff $d(x)=d(y)$.
The CCF Markov chain $p(k)$ and the DCF$^{(n)}$ Markov chain
$\ell(k)$ will be realized as
\begin{equation}\Label{rugby}
p(k):=\pi_c(c_k,d_k,e_k) \quad\mbox{and}\quad
\ell(k):=\pi_d(c_k,d_k,e_k).
\end{equation}\abel{rugby}
The flag $e_k$ indicates whether the coupling between the two
processes $p(k)$ and $\ell(k)$ is considered to be  still in force
($e=0$) or to have already broken down ($e=1$).

The distribution of $(c_0,d_0,e_0)$, that is the initial
distribution of the Markov chain, is defined as the image of $\mu$
under the function $\Phi^{(n)}=(\Phi_1^{(n)},\Phi_2^{(n)},0):
\Omega_1\longrightarrow\Omega_{cde}^{(n)}$ which assigns to each
element of $\Omega_1$  an equivalent function $c$ and an
approximating function $d$ as follows (see Figure \ref{uruguay}):
\begin{eqnarray*}
\Phi_1^{(n)}(p)(x)&:=&\sum_{j\geq 1}j \won
\left\{x\in[0,np_j)+n\sum_{i=1}^{j-1}p_i\right\}\hspace*{11mm}(x\in
[0,n))\\
\Phi_2^{(n)}(p)(m)&:=&\Phi_1^{(n)}(p)(m-1)\hspace*{40mm}(m\in\{1,\ldots,n\}).
\end{eqnarray*}
Thus $p(0)=p$ and $\ell(0)=\pi_d(\Phi^{(n)}(p))$ are the initial continuous
 and discrete partitions generated by  $p\in\Omega_1$.

 \begin{figure}[t]
 \begin{center}
 \epsfig{file=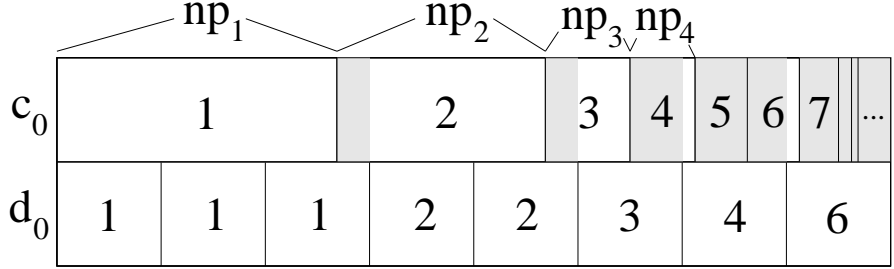,height=3.6cm,angle=0}
 \begin{minipage}[t]{12cm}
 \caption{ Constructing a continuous partition of $[0,n)$ and a
 discrete partition of $\{1,\ldots n\}$ from a partition
 $p\in\Omega_1$. Here $n=8$. The shaded area indicates the region
 where the continuous and the discrete numbering disagree.\Label{uruguay}}
 \end{minipage}\end{center}\vspace*{-5mm}
 \end{figure}\abel{uruguay}

\noindent
  {\bf Proof of (\ref{m8}):}
To bound $N_{\ell(0)}$, the number of parts in $\ell(0)$, observe
that all the pieces in $p$ of size less than $1/n$ can give rise
to at most $\sum_i np_i \won_{np_i<1}$ parts (singletons) in
$\ell(0)$. Therefore,
\begin{eqnarray}
N_{\ell(0)}&=&\#\Phi_2^{(n)}(p)[\{1,\ldots,n\}] \ =\ \nonumber
\sum_{i}\won_{np_i\geq 1}+ np_i \won_{np_i<1}\\
& \leq&
\sum_{i}(np_i)^\al\won_{np_i\geq 1}+ (np_i)^\al \won_{np_i<1}
\ =\
n^\al\sum_{i}p_i^\al.\Label{coseq}
\end{eqnarray}\abel{coseq}
Consequently, due to assumption (\ref{m7}),
\begin{equation}\Label{denmark}
E_\mu^{(n)}[N_{\ell(0)}]=O(n^\al)
\end{equation}\abel{denmark}
and hence, for all $\beta>\al$,
\[ Q_\mu^{(n)}\left[N_{\ell(0)}>n^\beta\right]\leq n^{-\beta}E_\mu^{(n)}[N_{\ell(0)}]
=O(n^{\al-\beta}) ,\]
 thus proving (\ref{m8}).

We now define informally the kernel of the Markov chain
$(c_k,d_k,e_k)$ with state space $\Omega_{cde}^{(n)}$. Assume that
the current state of the Markov chain is $(c,d,e)$. To compute the
state $(\bar{c},\bar{d},\bar{e})$ the Markov chain is going to
jump to in the next step we generate four random variables $\xi_1,
\xi_2, \zeta_1$ and $\zeta_2$ such that $\xi_1$ and $\xi_2$ and
$(\zeta_1,\zeta_2)$ are independent of each other and of
everything else and such that the $\xi_i$ are uniformly
distributed on $[0,n)$ and    $(\zeta_1,\zeta_2)$ is uniformly
distributed on $[0,n)^2\backslash \bigcup_{j=1}^n [j-1,j)^2$. The
$\xi_i$ will serve to sample uniformly with replacement from
$[0,n)$ whereas the $\zeta_i$ will be used to sample uniformly
without replacement from $\{1,\ldots,n\}$ in case the $\xi_i$ have
chosen the same atom in $d$ twice.
The new continuous partition $\bar{c}$  is then defined as follows:
\begin{eqnarray}
\mbox{If}\quad c(\xi_1)&\ne&c(\xi_2):\nonumber\\
\Label{cn}\bar{c}(x)&=&\left\{\begin{array}{ll} c(\xi_1) &\mbox{if
$c(x)=c(\xi_2)$}\\ c(x)&\mbox{else.}
\end{array}
\right.\\ \mbox{If}\quad c(\xi_1)&=&c(\xi_2):\nonumber\\
\Label{ce}\bar{c}(x)&=&\left\{\begin{array}{ll} {\rm new\,}(c,d)
&\mbox{if $c(x)=c(\xi_1)$ and $x> \xi_1$}\\ c(x)&\mbox{else.}
\end{array}
\right.
\end{eqnarray}
\abel{cn}\abel{ce}
We see that
  the two parts are indeed chosen with probabilities given by their size.
 In (\ref{cn}) two different sets, of sizes $\Leb(c^{-1}(c(\{\xi_i\})),\ i=1,2,$
 have been
selected and are merged by assigning the set $c^{-1}(c(\{\xi_2\})$,
hit by $\xi_2$, the
number $c(\{\xi_1\})$ of the set $c^{-1}(c(\{\xi_1\})$,  selected by $\xi_1$.
This creates a new set $\bar{c}^{-1}(c(\xi_1))$ with Lebesgue measure
$\Leb(c^{-1}(c(\xi_1)))+\Leb(c^{-1}(c(\xi_2)))$.\\
 In (\ref{ce}) the set
$c^{-1}(c(\{\xi_1\})=c^{-1}(c(\{\xi_2\})$ is
chosen twice, so it has to be split. Since $\xi_1$ is conditionally uniformly distributed
on this set we can reuse it as splitting point for that set:  The
part to the left of $\xi_1$ retains its old number $c(\xi_1)=c(\xi_2)$ whereas the
part to its right gets a new number ${\rm new\,}(c,d)$, which is
not in the range of $c$ or $d$. Note that it is always possible to
find such a new number since $\ZZ\backslash c[[0,n)]$ is assumed
to be infinite. By comparing this with the definition of CCF given at the
beginning of the Introduction we see that $p(k)$ defined in (\ref{rugby})
is a CCF Markov chain.

In the discrete case, the two parts chosen are the ones containing the numbers
$\lceil\xi_1\rceil$ and $\lceil\xi_2\rceil$, which ensures that the parts are
chosen size biased.
The rule for merges in the discrete partition is analogous to
(\ref{cn}):
\begin{eqnarray}
\mbox{If}\quad d(\lceil\xi_1\rceil)&\ne&
d(\lceil\xi_2\rceil):\nonumber\\
\label{dn}\bar{d}(m)&=&\left\{\begin{array}{ll}
d(\lceil\xi_1\rceil) &\mbox{if $d(m)=d(\lceil\xi_2\rceil)$}\\
d(m)&\mbox{else.}
\end{array}
\right.
\end{eqnarray}\abel{dn}
Here two different parts with numbers $d(\lceil\xi_1\rceil)$ and
$d(\lceil\xi_2\rceil)$ have been chosen. They are merged by giving both of them
the number  $d(\lceil\xi_1\rceil)$.\\
The rule for splitting is slightly more complicated. If the same
part (but not the same atom) is sampled twice by the $\xi_i$ then
again,  as in the continuous setting, $\xi_1$ determines the point at which the set $d^{-1}(\{\lceil
\xi_1\rceil\})$ is going to be split: The points to the left of
$\lceil \xi_1\rceil$
 and the points to the right of $\lceil \xi_1\rceil$
will constitute the two new fragments. The point $\lceil
\xi_1\rceil$ itself will be attached to the left or the right part
in such a way that the splitting rule  for DCF$^{(n)}$, given in (\ref{Kn}) and (\ref{Kn1}), is imitated.
This is done as follows:
\begin{eqnarray}
\mbox{If}\quad
d(\lceil\xi_1\rceil)&=&d(\lceil\xi_2\rceil)\quad\mbox{and}\quad
\lceil\xi_1\rceil\ne\lceil\xi_2\rceil:\nonumber\\ \Label{den}
\bar{d}(m)&=&\left\{\begin{array}{ll} {\rm new\,}(c,d) &\mbox{if
$d(m)=d(\lceil\xi_1\rceil)$ and $m>\lceil\xi_1\rceil$}\\
&\mbox{or}\  m=\lceil\xi_1\rceil\ \mbox{and} \
\xi_1<\lfloor\xi_1\rfloor+{\displaystyle
\frac{\#d^{-1}(\{d(\lceil\xi_1\rceil)\})\cap[0,\xi_1]}
{\#d^{-1}(\{d(\lceil\xi_1\rceil)\})-1}}\\ d(m)&\mbox{else}.\\
\end{array}\right.
\end{eqnarray}
Indeed, consider for simplicity the case that the atoms of
the set $d^{-1}(\{d(\lceil\xi_1\rceil)\})$ are not scattered around the whole set $\{1,\ldots,n\}$,
which they typically will be, but are collected  at the bottom:  $d^{-1}(\{d(\lceil\xi_1\rceil)\})=
\{1,\ldots,a\}$, where $a:=\sharp d^{-1}(\{d(\lceil\xi_1\rceil)\})$.  Definition (\ref{den}) tells us that this set is
split into $\{1,\ldots,j\}$ and $\{j+1,\ldots,a\}$ if
\[j-1+\frac{j-1}{a-1}\leq \xi_1< j +\frac{j}{a-1}.\]
Conditioned on $\xi_1\in[0,a)$, the probability for this to happen is $1/(a-1)$.
This means that the discrete set $\{1,\ldots,a\}$ is indeed split as described at the beginning
of Section \ref{DCF}.

If however  the same atom in
$d$ has been sampled twice by the $\xi_i$'s, i.e.\ $\lceil\xi_1\rceil=\lceil\xi_2\rceil$, then  $\xi_1$ and $\xi_2$ are
disregarded and  $\bar{d}(n)$ is defined as in (\ref{dn}) and
(\ref{den}) but with $(\xi_1,\xi_2)$ replaced by $(\zeta_1,\zeta_2)$
in order to sample without replacement.
The process $\ell(k)$ defined in
(\ref{rugby}) is a DCF$^{(n)}$ Markov chain.

It remains to define $\bar{e}$:
\[\bar{e}=\left\{\begin{array}{ll}
1&\mbox{if $\lceil \xi_1\rceil= \lceil \xi_2\rceil$ or
 $c(\xi_1)\ne d(\xi_1)$ or
 $c(\xi_2)\ne d(\xi_2)$}\\
e&\mbox{else.}
\end{array}
\right.\] In the case $\bar{e}=1$ the  coupling has broken down:
Either the same atom in the discrete partition has been sampled
twice by the $\xi_i$'s or at least one of the  $\xi_i$'s belongs to
non-corresponding sets in the continuous and the
discrete partition. The
 time
$\tau:=\inf\{k\geq 1:\ e_k=1\}$ is regarded as the decoupling time of the chains
$p(k)$ and $\ell(k)$.\\
The definition of the transition
kernel for the Markov chain on $\Omega_{cde}^{(n)}$ is now
complete. It is summarized in Figures \ref{mfig} to \ref{bfig}.

\begin{figure}[t]
\begin{center}
\epsfig{file=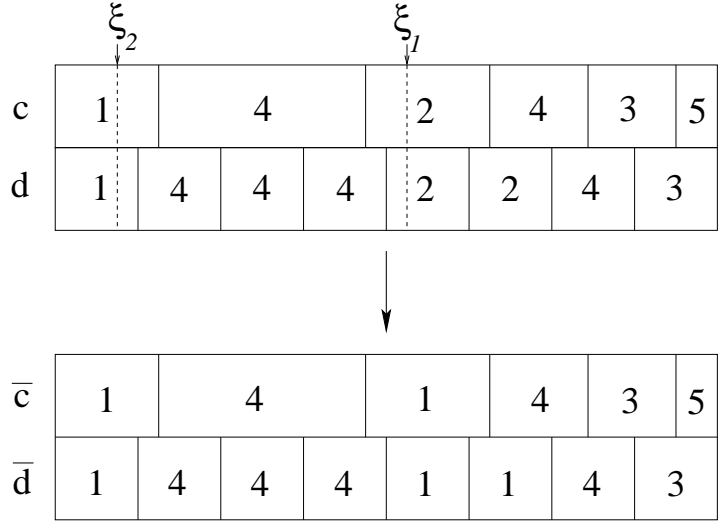,height=7cm,angle=0}
\caption{Merging the parts with numbers 1 and 2.\Label{mfig}}
\end{center}\vspace*{-5mm}
\end{figure}
\begin{figure}[h]

\begin{center}
\epsfig{file=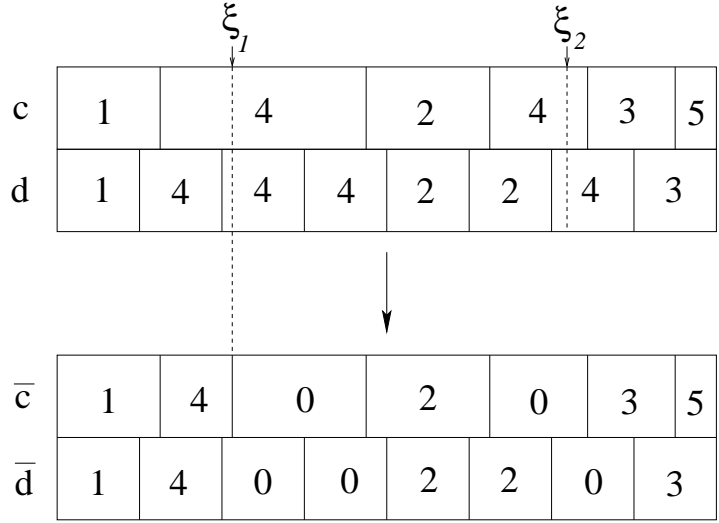,height=7cm,angle=0}
\begin{minipage}[t]{12cm}
\caption{Splitting the part with number 4 into a part with number
4 and a part with number ${\rm new}\,(c,d)=0$.\Label{sfig}}
\end{minipage}
\end{center}\vspace*{-5mm}
\end{figure}
\begin{figure}[h]
\begin{center}
\epsfig{file=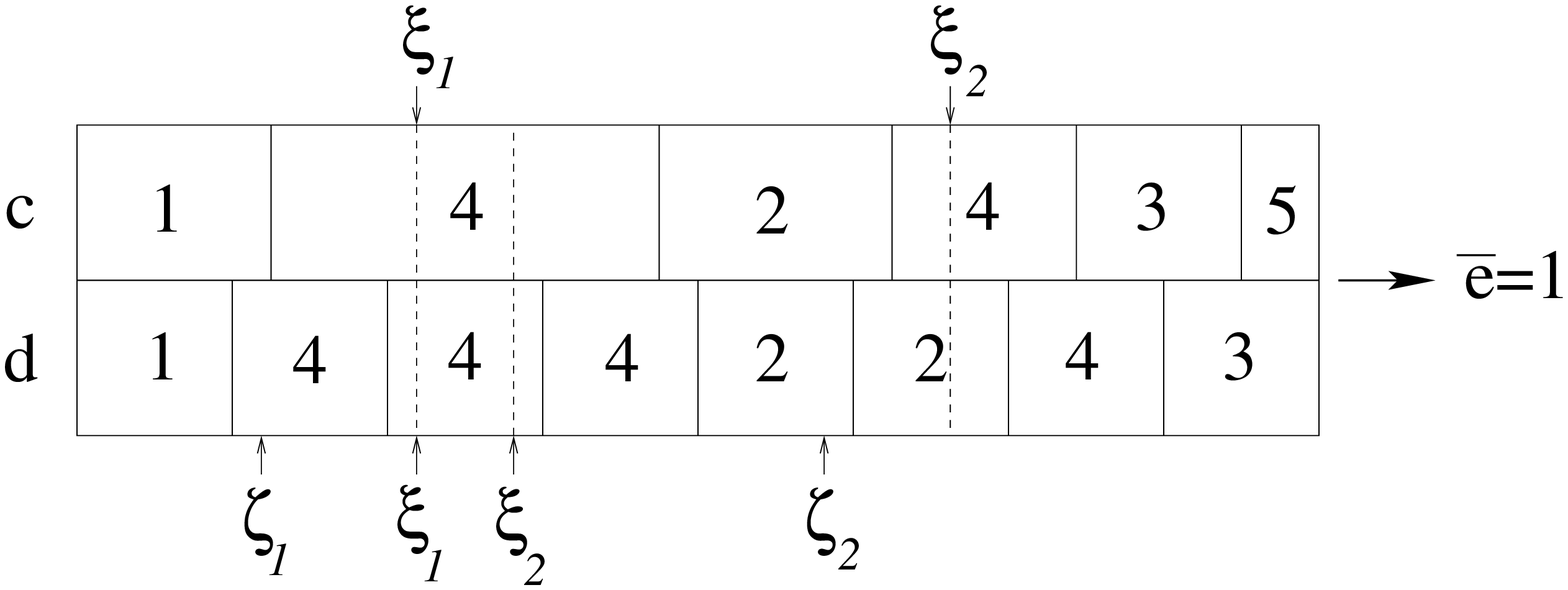,height=4.2cm,angle=0}
\vspace*{-5mm}
\begin{minipage}[t]{12cm}
\caption{Two ways to  decouple the chains: Sampling from the
 region where $c$ and $d$ disagree ($\xi_2$, top)
or sampling with replacement from $d$ ($\xi_1$ and $\xi_2$,
bottom).\Label{bfig}}
\end{minipage}
\end{center}
\end{figure}

\noindent
  {\bf Proof of (\ref{m9}):}
We denote by
\[\varrho_k:=\Leb(\{x\in[0,n)\ :\ c_k(x)\ne d_k(\lceil x\rceil\}).\]
the discrepancy between $c_k$ and $d_k$. For $k=0$, this
is the roundoff error caused by the approximation of $c_0$ by $d_0$; its
size is the length of the shaded area in Figure \ref{uruguay}.
Note that
\begin{equation}\Label{simple}
\varrho_0\leq N_{\ell(0)}
\end{equation}\abel{simple}
because any part in $d_0$ might disagree with $c_0$ at most in its
right most atom. Moreover, $\varrho_k$ can increase in each step
by at most 1 as long as $k<\tau$: Indeed, if two parts are merged,
$\varrho_k$ does not increase at all (it might even decrease)
whereas it might increase by at most $\Leb((\lfloor\xi_1\rfloor,
\lceil \xi_1\rceil])=1$ in case of splitting. Hence,
$\varrho_{k+1}\leq \varrho_{k}+1$ if $k<\tau$ and therefore,
\begin{equation}\Label{swiss}
  \varrho_k\leq\varrho_{0}+k\quad\mbox{on the event}\quad \{k<\tau\}.
\end{equation}\abel{swiss}
Since the $|\cdot|_1$-diameter of $\Omega_1$ is at most  2 we have
\begin{eqnarray}\nonumber
\lefteqn{E_\mu^{(n)}\left[\left|p(\lfloor n^\beta\rfloor)-
\frac{\ell(\lfloor n^\beta\rfloor)}{n}\right|_1\right]}\\
&\leq&\Label{nigeria} E_\mu^{(n)}\left[\left|p(\lfloor
n^\beta\rfloor)- \frac{\ell(\lfloor n^\beta\rfloor)}{n}\right|_1,
\lfloor n^\beta\rfloor \leq \tau \right]+ 2Q_\mu^{(n)}[\tau<\lfloor
n^\beta\rfloor].
\end{eqnarray}\abel{nigeria}
We are going to bound the first term in (\ref{nigeria}) first. It
is easy to see that $|p-q|_1\geq |{\rm sort\, }(p)-{\rm sort\,
}(q)|_1$ for any two summable sequences $p=(p_i)_i$ and
$q=(q_i)_i$ of non-negative numbers. Indeed, if $p_i>p_j$ and
$q_i<q_j$, then swapping $q_i$ and $q_j$ would not increase $|p-q|_1$.
Therefore, by definitions (\ref{france}), (\ref{senegal}) and (\ref{rugby}) on
the event $\{\lfloor n^\beta\rfloor<\tau\}$,
\begin{eqnarray*}\left|p(\lfloor n^\beta\rfloor)-
\frac{\ell(\lfloor n^\beta\rfloor)}{n}\right|_1 &\leq&\frac{1}{n}\sum_{i\geq
1}\left|\Leb (c_{\lfloor n^\beta\rfloor}^{-1}(\{i\}))-\#
d_{\lfloor n^\beta\rfloor}^{-1}(\{i\}) \right|\\ &\leq
&\frac{1}{n}\sum_{i\geq 1}\Leb( \{x:
i\in\{c_{\lfloor n^\beta\rfloor}(x),d_{\lfloor n^\beta\rfloor}(x)\},
c_{\lfloor n^\beta\rfloor}(x)\ne d_{\lfloor n^\beta\rfloor}(x)\})\\ &\leq&
\frac{2}{n}\varrho_{\lfloor n^\beta\rfloor}\ \leq\
\frac{2}{n}(\varrho_0+\lfloor n^\beta\rfloor)\ \leq\ \frac{2}{n}(N_{\ell(0)}+\lfloor n^\beta\rfloor)
\end{eqnarray*}
by (\ref{swiss}) and (\ref{simple}).
Consequently, due to (\ref{denmark}), the first term in
(\ref{nigeria}) is of order $O(n^{\al-1}+n^{\beta-1})$, thus going to
0 as $n\to\infty$.

To show that the  second term in (\ref{nigeria}) goes to 0 as well
we assume without loss of generality that $\al<\beta<1/2$.
Consider
the probability that a chain which has not decoupled until the
$k$th step will decouple in the $(k+1)$th step. Given $\varrho_0,\ldots,\varrho_k$,
the event that $\xi_1$ samples two different parts in $c_k$ and
$d_k$ has probability $\varrho_k/n$. The same holds for $\xi_2$.
Moreover, the event that one atom in $d_k$ is sampled twice,
i.e.\ that  $\lceil \xi_1\rceil=\lceil \xi_2\rceil$ has probability $1/n$.
Therefore, the probability that either of these events occurs and
the chain decouples is at most $(2\varrho_k+1)/n$. On the event
$\{\tau>k, \varrho_0< n^\beta\}$ this can be bounded from above due
to (\ref{swiss}) by $(2(n^\beta+k)+1)/n$ which is less than
$5n^{\beta-1}$ if $k\leq n^\beta$. Thus we get by induction over
$k$,
\[Q_\mu^{(n)}[\tau>k, \varrho_0< n^\beta]\geq (1-5n^{\beta-1})^k
Q_\mu^{(n)}[\varrho_0< n^\beta]\] for all $k\leq n^\beta$ and hence
 \begin{eqnarray}
Q_\mu^{(n)}[\tau\geq \lfloor n^\beta\rfloor]&\geq&
\left(\left(1-5n^{\beta-1}\right)^{n^{1-\beta}}\right)^{n^{2\beta-1}}
Q_\mu^{(n)}\left[\varrho_0< n^{\beta}\right].\Label{bren}
\end{eqnarray}\abel{bren}
Due to $2\beta-1<0$, the first factor in (\ref{bren}) converges to
one as $n\to\infty$. The same holds for the second factor due to
(\ref{simple}) and (\ref{m8}). Consequently, also the second term
in (\ref{nigeria}) goes to 0, which completes the proof of
(\ref{m9}).\qed
\end{proof}

\section{\DCFn\ convergence}   \Label{DCFconv}
 It was mentioned in the Introduction, that the {\it uniform} rate of
 convergence to $\pins$ is too weak to combine properly with $n\to\infty$.
 However, according to the following theorem (to be proved in
 subsection~\ref{proofofer}), the situation is better when starting off from
 partitions with relatively few parts and restricting our attention to a certain
 family $\calC$ of $\Omega_1$--neighborhoods to be defined below. For every
 $n\in\Nat$ and $\beta\in(0,1]$, thus, denote accordingly
 \begin{eqnarray*}
   \Pnb&=&\left\{\ell\in\Pn\,:\ N_\ell<n^\beta\right\}
        = \left\{\ell\in\Pn\,:\ \ell_{\lceil n^\beta\rceil}=0\ \right\}\ .
 \end{eqnarray*}
 As for the definition of $\calC$, for each $k\!\in\!\Nat$ let
 \begin{equation}  \Label{Ik}
  I_k=\left\{(\aaa,\bbb)\!=\!(a_i,b_i)_{i=1}^k\,:\ \ \
0<a_i<b_i<1\,\ \ \
      \sum_{i=1}^kb_i\!<\!1,\ \ \ \ \ a_k\!>\!1\!-\!\sum_{i=1}^ka_i\ \right\}
 \end{equation}
 and denote $\dab\!=\!\min\left\{1\!-\!\sum_{i=1}^kb_i
             \ \ ,\ \ a_k\!-\!(1\!-\!\sum_{i=1}^ka_i)\right\}$.
 Then, for each $(\aaa,\bbb)\in I_k$, define
 \begin{equation}  \Label{Cab}
   \Cab=\left\{x\in\Omega_1\,:\ {}_{_{_{}}}^{^{^{}}}
           x_i\in (a_i,b_i)\ \ \ \mbox{for}\ i=1,\ldots,k\right\}\hspace{5cm}
 \end{equation}
 which is nonempty if and only if\ \ \
 $0\!<\!a_i\!\!<\!\!\underset{1\le j\le i}{\min}b_j$\ \ \ for $i=1,\ldots,k$,
 \ \ in which case the conditions on $(\aaa,\bbb)\!\in\!I_k$ guarantee that
 \begin{equation}\Label{altC}
  \Cab\!=\!\left\{x=(x',x'')\ :\ x'\in G_{{\saa,\sbb}_{_{}}}^{^{}},
         \ x''\in (1-|x'|_1)\Omega_1\right\}.\hspace{4cm}
 \end{equation}
 (Here $(\,\cdot\,,\,\cdot\,)$ denotes concatenation and $G_{\saa,\sbb}$ is the
 (nonempty) subset of points in $\prod_{_{i=1}}^{^k}\!(a_i,b_i)$ whose coordinates
  are nonincreasing). Moreover, $I_k$'s definition~(\ref{Ik}) implies that
 \begin{equation}\Label{x'x''cond}
  \dab<|x''|_1<x'_k-\dab\ ,\hspace{1.5cm}\forall (x',x'')\in\Cab\ .\hspace{4.5cm}
 \end{equation}
 Finally,
\begin{equation} \Label{calC}
 \calC=\{\Cab\, :\, (\aaa,\bbb)\in I_k,\ k\geq 1\}.
 \end{equation}
 The family $\calC$ of $\Omega_1$--neighborhoods will be shown in
 Section~\ref{Vershik} to be sufficiently rich to characterize
 $\wmu_1$ uniquely. At the same time, and as a result of their
 special features~(\ref{altC}) and~(\ref{x'x''cond}), the convergence
 of the \DCFn\ to its equilibrium is fast on the sets in $\calC$:
 \begin{theorem} \Label{the-xxx}
 Fix $\beta\!\in\!(0,\frac{1}{2})$. For each $n\!\in\!\Nat$ let
 $\left(X^{(n)}(k)\right)_{k\ge 0}$ be a DCF${}^{(n)}$ Markov
 chain with underlying probability measure $P^{(n)}$ and initial
 distribution $\mu_0^{(n)}\!\in\!\calM_1(\calP_{n,\beta})$.
 Then for any $C\!\in\!\calC$,\ $\beta'\!>\!\beta$ and integer sequence
 $k=k_n\ge n^{\beta'}$
 \[ \dnc:=P^{(n)}\left(X^{(n)}(k)\in nC\right)-\pins(nC)
                \underset{n\to\infty}{\longrightarrow} 0\,.\]
 \end{theorem}
  \subsection{Characters in $S_n$ -- Background}    \Label{characters}
  \noindent
  Recall that the partition space $\Pn$ can be viewed as the quotient of the
  permutation group $S_n$ under conjugacy. Thus the natural inner product on
  $F_n\!:=\!\{f:\Pn\longrightarrow\Reals\}$\ is
  \[ \lip f,g\rip=\lip  f,g\rip_{_n}=\sum_{\ga\in\Pn}f(\ga)g(\ga)\pins(\ga)\ .\]
  The fact mentioned earlier that $\pins$ is a reversing measure for the \DCFn\
  means precisely that $\Kn$ is selfadjoint with respect to this inner product.
  \medskip

  \noindent
  The following basic facts regarding the character theory of $S_n$, as well as the
  full theory, can be found, for example, in~\cite{jameskerber}, and their relevance
  to random group actions (such as transpositions in our case) in~\cite{diaconis} and
  ~\cite{flattoetal}.
  The characters $\{\chi\}$ of $S_n$ (traces of the irreducible representations)
  are functions on $S_n$, constant on conjugacy classes, and as such can be seen to
  be functions on $\Pn$. They are orthonormal under $\lip\,\cdot\,,\,\cdot\,\rip$ and
  since there are $\#\Pn$ of them, they are indexed by the partitions
  (\,$(\Xla)_{_{\la\in\Pn}}$) and form an orthonormal base of $F_n$.
  \smallskip

  \noindent
  Since $\Kn$ represents a random transposition, its dual $K^{(n)^*}$ acts on
  $\calM_1(S_n)$ as a convolution
  \[ K^{(n)^*}\mu=\kappa^{(n)}\star\mu\hspace{2cm}
      \left(\kappa^{(n)}(\mbox{\small{transposition}})\!=\!\frac{2}{n(n\!-\!1)}\ \
                              \ \mbox{and}\ \ 0\ \ \mbox{otherwise}\right) \]
  as a result of which, and of a corollary (\cite[Ch. 2, Prop.~6]{diaconis})
  of Schur's lemma,
  \begin{itemize}
   \item[\rm \bf a)] $\Kn$'s eigenfunctions are the characters $(\Xla)_{\la\in\Pn}$
    \vspace{-.35cm}

   \item[\rm \bf b)] the eigenvalue $\thelan$ corresponding to $\Xla$is given by
        $\frac{\mbox{\large{$\Xla$}}(\mbox{\small{transposition}})}
              {\mbox{\large{$\Xla$}}(\mbox{\small{identity}})}$\ .
  \end{itemize}
  A result of Frobenius in principle provides formulae for all characters. Although
  in general they can be intractable, this is not so at transpositions and at the
  identity, thus yielding (\cite[D-2,p.40]{diaconis})
 \begin{equation} \Label{4.1}
  \thelan=\frac{1}{n(n-1)} \sum_j \la_j (\la_j - 2j+1)
         =\frac{1}{n(n-1)}\left(\sum_{i=1}^n \la_i^2-\sum_{j=1}^{\la_1}{\la_j'}^2
                                                                           \right)\,.
 \end{equation}
 ($\la$'s adjoint partition $\la'$ is defined below). In particular
 $\theta^{(n)}_{(n,0,\ldots)}=1$ and $\chi_{_{(n,0,\ldots)}}\equiv 1$.
 \bigskip

 \noindent
 For many purposes, a partition $\la\in\Pn$ can be best described by its
 Young diagram $\Ul$ (Fig.~\ref{Mnfig}), consisting of $N_\la$ rows of
 $\la_1,\ldots,\la_{N_\la}$ cells respectively, in terms of which some
 useful features of $\la$ can be defined. The $j$-th cell in row $i$
 is denoted $(i,j)$.
 \begin{itemize}
  \item $\la'\in\Pn$ is the partition whose Young diagram is obtained from
        $\la$'s by transposition; $\Ups_{_{\!\la'}}=\Ul^{^{\mbox{\tiny $T$}}}$
  \item $B_\la=\max\{i:(i,i)\in\Ul\} =\max\{i:\la_i\ge i\}$
            \hspace{1cm}($\la$'s diagonal length)
  \item $R_\la(i,j)=\{(u,v)\,:\,\ i\!\le\!u\!\le\!\lambda_j',\ \
                 j\,\mbox{\footnotesize{$\vee$}}\,\la_{u+1}\!\le\!v\!\le\!\la_u\}$
            \hspace{1cm}($\Ul$'s rim segment straddled by $(i,j)$)
  \item $\Ups_{\!\lij}=\Ul\setminus  R_\la(i,j)$ defines $\lij$ \hspace{.3cm}
   (\begin{minipage}[t]{9.3cm}
       a diagram obtained from $\la$'s by removing a rim segment is a Young
       diagram; this defines the partition $\lij$)
   \end{minipage}
 \end{itemize}

 \noindent
 In addition, for any $\ga\in\Pn$, define
 $\ga^{\widehat{r}}=(\ga_1,\ldots,\widehat{\ga_r},\ldots)\in\calP_{n-\ga_r}$,
 the partition obtained from $\ga$ by removing its $r$-th part.
 \noindent
 The following \MN\ rule (see \cite[Theorem 3.4]{flattoetal}) provides a way of
 recursively evaluating characters: for all $\la,\ga\in\Pn$\  and\
 $1\!\le\!r\!\le\!N_\ga$
 \begin{equation}  \Label{MN}
  \Xla(\ga)  = \hspace{-.4cm}\sum_{(i,j)\,:\,\#R_\la(i,j)=\ga_r}\hspace{-.4cm}
         (-1)^{\la'_j-i}\,\chi_{_{_{\!\lij}}}\!(\ga^{\widehat{r}})\hspace{5cm}
 \end{equation}
 in the sense that the sum is zero if its index set is empty, and
 $\chi_{_{\emptyset}}(\mbox{\footnotesize{$\emptyset$}})=1$. Thus, for a fixed
 order in which $\ga$'s parts are chosen, $\Xla(\ga)$ can be calculated by
 covering all possible ways of successively stripping off $\ga_r$--sized rim
 segments from $\la$'s diagram, and $\Xla(\ga)=0$ if it is impossible to
 exhaust $\Ul$ entirely in this way. In particular
 \begin{equation} \Label{NvsB}
  N_\ga<B_\la\Longrightarrow \Xla(\ga)=0\hspace{7cm}
 \end{equation}
 since any rim segment of $\la$ contains at most one diagonal cell $(i,i)$.

   \begin{figure}
 \begin{picture}(450,160)(0,-20)
  \put(5,5){\line(1,0){15}}
  \put(5,20){\line(1,0){15}}
  \put(5,35){\line(1,0){60}}
  \put(5,50){\line(1,0){60}}
  \put(5,65){\line(1,0){105}}
  \put(5,80){\line(1,0){120}}
  \put(5,95){\line(1,0){120}}
  \put(5,110){\line(1,0){120}}
  \put(5,110){\line(0,-1){105}}
  \put(20,110){\line(0,-1){105}}
  \put(35,110){\line(0,-1){75}}
  \put(50,110){\line(0,-1){75}}
  \put(65,110){\line(0,-1){75}}
  \put(80,110){\line(0,-1){30}}
  \put(95,110){\line(0,-1){30}}
  \put(110,110){\line(0,-1){45}}
  \put(125,110){\line(0,-1){30}}
%
  \thicklines
  \put(94.8,95.8){\line(1,0){30.5}}
  \put(50,80.8){\line(1,0){45}}
  \put(110,79){\line(1,0){15}}
  \put(65,64){\line(1,0){45}}
  \put(35,50.8){\line(1,0){15}}
  \put(35,34){\line(1,0){30.4}}
  \put(94,96.2){\line(0,-1){16.2}}
  \put(125.6,96){\line(0,-1){17}}
  \put(49,81){\line(0,-1){31}}
  \put(110.7,80){\line(0,-1){16}}
  \put(66,65){\line(0,-1){31}}
  \put(34.5,51.3){\line(0,-1){17.5}}
 \thinlines
 \put(5,110){\line(1,-1){60}}
  \put(43,88){\circle*{3}}
  \put(58,103){\circle*{3}}
  \multiput(48,88)(7,0){11}{\line(1,0){3}}
  \multiput(43,83)(0,-6.9){7}{\line(0,-1){3}}
 \put(1,160){\mbox{\small{$\la\ =(8,8,7,4,4,1,1,0,\ldots)$}}}
 \put(1,145){\mbox{\small{$\la'=(7,5,5,5,3,3,3,2,0,\ldots)$}}}
 \put(20,120){\mbox{\small{$B_\la\!=\!4$}}}
 \put(34,116){\vector(-1,-4){6.7}}
 \put(70,69){\mbox{\small{$R_\la(2,3)$}}}
 \put(80,48){\mbox{\small{$\lambda_*^{(2,3)}=(8,6,3,3,2,1,1,0,\ldots)$}}}
 \put(80,34){\mbox{\small{$\lambda_*^{(1,4)}=(7,6,3,3,3,1,1,0,\ldots)$}}}
  \put(265,20){\line(1,0){30}}
  \put(265,35){\line(1,0){30}}
  \put(265,50){\line(1,0){45}}
  \put(265,65){\line(1,0){105}}
  \put(265,80){\line(1,0){135}}
  \put(265,95){\line(1,0){150}}
  \put(265,110){\line(1,0){150}}
  \put(265,110){\line(0,-1){90}}
  \put(280,110){\line(0,-1){90}}
  \put(295,110){\line(0,-1){90}}
  \put(310,110){\line(0,-1){60}}
  \put(325,110){\line(0,-1){45}}
  \put(340,110){\line(0,-1){45}}
  \put(355,110){\line(0,-1){45}}
  \put(370,110){\line(0,-1){45}}
  \put(385,110){\line(0,-1){30}}
  \put(400,110){\line(0,-1){30}}
  \put(415,110){\line(0,-1){15}}
  \multiput(265,95)(15,0){9}{\line(1,-1){15}}
  \multiput(265,80)(15,0){9}{\line(1,1){15}}
 \put(261,160){\mbox{\small{$\ga\ =(10,9,7,3,2,2,0,\ldots)$}}}
 \put(261,145){\mbox{\small{$\ga'=(6,6,4,3,3,3,3,2,2,1,0,\ldots)$}}}
 \put(350,50){\mbox{\small{$\ga^{\widehat{2}}=(10,7,3,2,2,0,\ldots)$}}}
 \put(68,-15){\MN\ rule:
         $\Xla(\ga)= \chi_{_{\la_*^{(1,4)}}}(\ga^{\widehat{2}})
                    -\chi_{_{\la_*^{(2,3)}}}(\ga^{\widehat{2}})$}
 \end{picture}
 \hspace*{1.5cm}
 \begin{minipage}[t]{12cm}
      \caption{\Label{Mnfig}Young diagrams of $\la,\ga\in\calP_{33}$.
        Two $\la$-cells, $(1,4)$ and $(2,3)$, generate rim
        segments of size $9$, the latter shown explicitly,
        which the \MN\ rule ``peels off" together with the
        deletion of $\ga_2$.}
 \end{minipage}
 \end{figure}
 \subsection{Proof of Theorem~\ref{the-xxx}}   \Label{proofofer}
 \noindent
 Before proceeding with the proof itself, it will be helpful to characterize
 the $\ga\in\Pn$ which belong to $nC\!=\!nC_{\saa,\sbb}$ for given
 $k\!\in\!\Nat$ and $(\aaa,\bbb)\!\in\!I_k$ (assuming
 $C\mbox{\small{$\ne\emptyset$}}$). It follows from $\Cab$'s
 description~(\ref{altC}) that any such $\ga$ can be expressed as a
 concatenation $(\ga',\ga'')$ where $\ga'\!\in\!\Gabn$ and
 $\ga''\!\in\!\calP_{n-|\ga'|_1}$, and where $\Gabn$ consists of
 nonincreasing integer valued $k$-sequences $\ga'$ which by virtue
 of~(\ref{x'x''cond}) satisfy
 \begin{equation} \Label{ab}
  \mbox{\bf i)}\  |\ga'|_1<n \hspace{1.5cm}
  \mbox{\bf ii)}\ \exists \delta\!=\!\delta(C)\!>\!0\ \ \mbox{such that}\ \
                                    \ga'_k>(n\!-\!|\ga'|_1)+ \delta n\ .
 \end{equation}
 This state of affairs is illustrated in Figure~\ref{nCab}.

 \paragraph{Proof of Theorem~\ref{the-xxx}}
 Fix $C\in\calC$ and define $f_n=\won_{nC}$.
 Then, in terms of $\mu_0^{(n)}$'s density
    $\gn(\ga)=\frac{\mu_0^{(n)}(\ga)}{\pins(\ga)}$:
 \begin{eqnarray*}
  P^{(n)}\left(X^{(n)}(k)\!\in\!nC\right)
   \!\!&=&\!\!\sum_\ga\mu_0^{(n)}(\ga)\Knk f_n(\ga)=\lip\gn,\Knk f_n\rip
     =\sum_{\la\in\Pn}\!\!
                \thelank\lip\gn,\Xla\rip\lip f_n,\Xla\rip\\
  \mbox{and,\ \ since}\ \theta^{(n)}_{(n,0,\ldots)}=1
                            &\mbox{and}&\chi_{(n,0,\ldots)}\equiv 1,\\ \\
   \pins(nC)&=&\lip f_n,1\rip
     =\theta_{(n,0,\ldots)}^{(n)^k}\lip\gn,\chi_{_{(n,0,\ldots)}}\rip\lip
                     f_n,\chi_{_{(n,0,\ldots)}}\rip
  \end{eqnarray*}
  so that
 \begin{equation}  \Label{Deltasum1}
    \dnc=\sum_{(n,0,\ldots)\ne\la\in\Pn}
      \thelank\lip\gn,\Xla\rip\lip f_n,\Xla\rip\ .\hspace{4.1cm}
 \end{equation}
 By assumption, $\gn(\ga)=0$ whenever $N_\ga>n^\beta$. On the other hand,
 $\Xla(\ga)=0$ whenever $\Bla>n^\beta$ and $N_\ga\le n^\beta$ by the
 consequence~(\ref{NvsB}) of \MN's rule.
 Thus~(\ref{Deltasum1}) becomes
 \begin{equation*}
    \dnc=\!\!\sum_{\underset{B_\la\le
    n^\beta}{(n,0,\ldots)\ne\la\in\Pn}}\!\!
                        \thelank\lip\gn,\Xla\rip\lip f_n,\Xla\rip\ .
 \end{equation*}
   \begin{figure}
 \begin{picture}(450,80)(-8,0)  \Label{gatyp}
 \setlength{\unitlength}{2pt}
 \put(140,72){\mbox{\LARGE{$\gamma=(\gamma',\gamma'')$}}}
 \thicklines
 \put(0,80){\line(1,0){120}}
 \put(120,80){\line(0,-1){7}}
 \put(120,73){\line(-1,0){8}}
 \put(112,73){\line(0,-1){7}}
 \put(112,66){\line(-1,0){2}}
 \put(110,66){\line(0,-1){7}}
 \put(110,59){\line(-1,0){9}}
 \put(101,59){\line(0,-1){14}}
 \put(101,45){\line(-1,0){11}}
 \put(90,45){\line(0,-1){7}}
 \put(90,38){\line(-1,0){3}}
 \put(87,38){\line(0,-1){7}}
 \put(87,31){\line(-1,0){87}}
 \put(0,31){\line(0,1){49}}
 \put(-7,76){$\gamma'_1$}
 \thinlines
 \put(-7,34){$\gamma'_k$}
 \put(0,38){\line(1,0){87}}
 \multiput(0,31)(5,0){13}{\line(1,1){7}}
 \put(68,31){\line(0,1){7}}
 \put(50,55){$\ga'$}
 \put(0,31){\line(0,-1){31}}
 \put(0,0){\line(1,3){3}}
 \put(3,9){\line(2,3){4}}
 \put(7,15){\line(1,1){5}}
 \put(12,20){\line(2,1){8}}
 \put(20,24){\line(2,3){5}}
 \put(0,5){\line(1,0){2}}
 \put(0,8){\line(1,0){3}}
 \put(0,11){\line(1,0){4.5}}
 \put(0,14){\line(1,0){6.5}}
 \put(0,17){\line(1,0){9.2}}
 \put(0,20){\line(1,0){12.5}}
 \put(0,23){\line(1,0){18}}
 \put(0,26){\line(1,0){21}}
 \put(0,29){\line(1,0){23}}
 \put(18,12){$\ga''\ne$}
 \put(34,13.5){\circle{2.5}}
 \put(32.7,11){\line(1,2){2.5}}
 \put(17,16){\vector(-1,1){8}}
 \put(70,12){\fbox{\large{$\gamma'_k>|\gamma''|_1$}}}
 \put(65,15){\vector(-3,1){20}}
\end{picture}
 \hspace*{1.5cm}
 \begin{minipage}[t]{12cm}
 \caption{\Label{nCab} A partition $\ga$ in $n\Cab$ splits into
          its first $k$ rows $\ga'$ and the remainder $\ga''$ which is
          nonempty but smaller in size than $\ga'$'s last row.}
 \end{minipage}
 \end{figure}
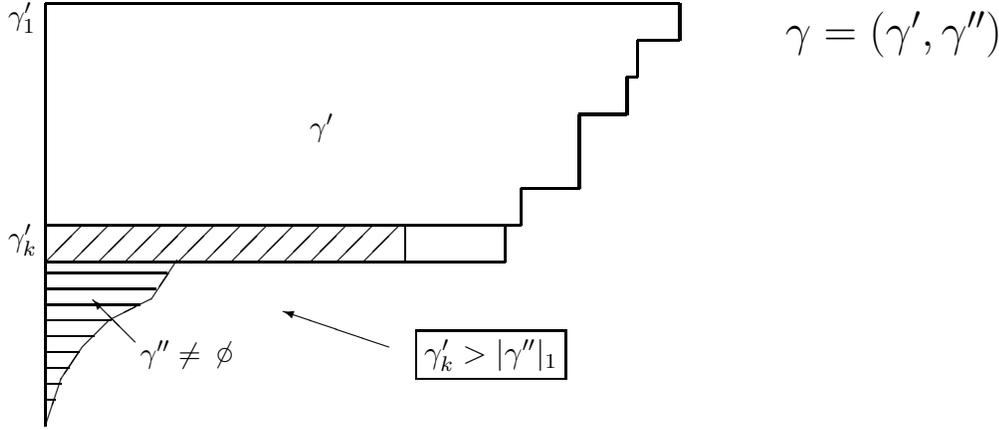

 Now choose an $\eta$ such that $1\!-\!(\beta'\!\!-\!\!\beta)<\eta<1$ and
 let $n_0=5^{\frac{1}{1-\eta}}$. Then, for all $n\!\ge\!n_0$,
 \begin{equation}  \Label{Deltasum2}
    \dnc=\left(\sum_{\la\in\Pn'}+\sum_{\la\in\Pn''}+\sum_{\la\in\Pn'''}\right)
                        \thelank\lip\gn,\Xla\rip\lip f_n,\Xla\rip
 \end{equation}
 where 
 \begin{eqnarray*}
  \Pn'  =\Pn'(\eta,\beta)  &=&\left\{\la\in\Pn:\ \ B_\la\le n^\beta,\ \ \
            \la_1\!\le\!n\!-\!2n^\eta,\ \ \ N_\la\!\le\!n\!-\!2n^\eta\,\right\} \\
  \Pn'' =\Pn''(\eta,\beta) &=&\left\{\la\in\Pn:\ \ B_\la\le n^\beta,\ \ \
                                     n\!-\!2n^\eta\!<\!\la_1\!<\!n\,\right\} \\
  \Pn'''=\Pn'''(\eta,\beta)&=&\left\{\la\in\Pn:\ \ B_\la\le n^\beta,\ \ \
                n\!-\!2n^\eta\!<\!N_\la\,\right\}
 \end{eqnarray*}
 (Our choice of $n_0$ ensures that $\Pn''$ and $\Pn'''$ are disjoint and that
  $(n,0,\ldots)\not\in\Pn'''$).\ It turns out that for $n$ large enough, the
  terms in~(\ref{Deltasum2}) vanish for all $\la\in\Pn''\cup\Pn'''$, whereas when
  $\la\in\Pn'$ the factor\ \ $|\thelan|$\ \ is sufficiently separated from $1$\ :
\begin{lemma} \Label{zeroterms} $\exists n_1=n_1(\beta,C)$ such that\ \
 $\lip f_n,\Xla\rip=0,\ \ \ \ \forall n\!\ge\!n_1,\
                              \forall \la\!\in\!\Pn''\cup\Pn'''$\ .
 \end{lemma}
 \begin{lemma} \Label{eigenvalues}
  For all $\la\!\in\!\Pn$,\ \
  $|\thelan|\le \bfrac{\la_1\mbox{\footnotesize{$\vee$}}N_\la}{n}$
  \ \ and thus\ \  $\exists n_2=n_2(\eta)$ such that
  \begin{equation}  \Label{thetabound}
   |\thelan| \le e^{-n^{\eta-1}}\hspace{1.5cm}\forall n\!\ge\!n_2,
                                       \ \ \forall\la\!\in\!\Pn'  \ .
  \end{equation}
 \end{lemma}
 \paragraph{Proof of Lemma~\ref{zeroterms}}
 Consider first $\la\!\in\!\Pn''$. Now, $C\!=\!\Cab$ for some
 $k\!\in\!\Nat$ and $(\aaa,\bbb)\!\in\!I_k$, so that, as discussed at the beginning of the
 section and illustrated in Figure~\ref{nCab}, $\ga$ can be split into
 $(\ga',\ga'')$ and
 \begin{equation} \Label{doublesum}
   \lip f_n,\Xla\rip=\sum_{\ga'\in\Gabn}
               \,\sum_{\ \ga''\in\calP_{n-|\gamma'|_1}}
               \!\!\!\pins(\ga',\ga'')\Xla(\ga',\ga'')\ .
 \end{equation}
 (Note that property~(\ref{ab}{\bf i}) guarantees that the inner sum is not
 vacuous,  i.e. $|\ga''|_1>0$).
 \smallskip

 \noindent
 We shall show that for every fixed $\ga'\!\in\!\Gabn$ the inner sum
 in~(\ref{doublesum}) equals zero. First apply \MN's rule~(\ref{MN}) $k$ times
 to $\Xla(\ga',\ga'')$ by successively stripping rim segments from $\la$, of
 lengths $\ga'_i$ at each stage $i,\ i=1,\ldots,k$. On the one hand
 $\lambda_1\!>\!(1\!-\!\delta)n$ for $n\ge n_1''=n_1''(\beta,\eta,C)$
 (since $\la\!\in\!\Pn''$), and on the other $\ga'_i>\delta n,\ i=1,\ldots,k$
 (by~(\ref{ab}{\bf ii})). This implies that at each of these $k$ reduction
 stages precisely one rim segment can be stripped off, namely the last $\ga'_i$
 cells of whatever remains of $\la_1, i=1,\ldots,k$. Summing up
 \begin{equation}  \Label{Mnreduction}
  \Xla(\ga',\ga'')= \Xlast(\ga'')
 \end{equation}
  where $\lambda^*\in\calP_{n-|\ga'|_1}$ is defined by
  $\lambda^*_1\!=\!\la_1\!-\!|\ga'|_1$ and $\lambda^*_j\!=\!\la_j,\ j\ge 2$.
  As for the first factor of the summand in~(\ref{doublesum}), note
  that~(\ref{ab}{\bf ii}) implies $\ga'_k>\ga''_1$ (see Figure~\ref{nCab})
  and thus
 \begin{equation} \Label{stationary}
  \frac{\pins(\ga',\ga'')}{\pi_S^{(n-|\ga'|_1)}(\gamma'')} =
  \frac{1}{\prod_{i=1}^k \ga'_i \prod_j N_{\ga'}(j) !} =:R(\ga').
 \end{equation}
 Inserting~(\ref{Mnreduction}) and~(\ref{stationary}) in the inner sum
 of~(\ref{doublesum}) we obtain
 \[ \sum_{\ga''\in\calP_{n-|\gamma'|_1}}\pins(\ga',\ga'')\Xla(\ga',\ga'')
      =R(\ga')\lip\Xlast,1\rip_{_{n-|\ga'|_1}}=0 \]
 since $\lambda^*$ is not the trivial partition, that is
 $\lambda^*\ne (n-|\ga'|_1,0,\ldots)$, (because $\lambda\ne (n,0,\ldots)$),
 and thus $\Xlast$ is orthogonal to $\chi_{(n-|\ga'|_1,0,\ldots)}\equiv 1$.

 \noindent
 The proof for $\la\in \Pn'''$ is similar, with
$n\!\ge\!n_1'''=n_1'''(\beta,\eta,C)$,
 where now the only rim segments which can be stripped off from $\la$ are from
 its first column. It remains to define $n_1\!=n_1''\mbox{\small{$\vee$}}n_1'''$.
 \qed
 \paragraph{Proof of Lemma~\ref{eigenvalues}}
 Using the formula for $\thelan$ given in~(\ref{4.1}),
 \begin{eqnarray*}\ \Label{4.1calc}
  \thelan &=& \frac{1}{n(n\!-\!1)}\,\left(\sum_{i=1}^n\la_i^2
                                               -\sum_{j=1}^{\la_1}{\la_j'}^2\right)
        \le  \frac{1}{n(n\!-\!1)}\,(\lambda_1 n -\lambda_1)=\frac{\lambda_1}{n}\ ,
 \end{eqnarray*}
 whereas, by duality, \
 $-\thelan=\theta_{\la'}^{(n)}\le \bfrac{\la'_1}{n}=\bfrac{N_\la}{n}$.\ \
 Moreover, for $\la\in\Pn'$,\ \ \
 \begin{eqnarray*}
 |\thelan|&\le&\left(1-\frac{2}{n^{1-\eta}}\right)
  =\left(1-\frac{2}{n^{1-\eta}}\right)^{n^{1-\eta}\ n^{\eta-1}}\!\!\!
  \le e^{-n^{\eta-1}}\ \ \ \ \ \mbox{as soon as}\ \
  \left(1\!-\!\frac{2}{n^{1-\eta}}\right)^{n^{1-\eta}}\!\!\le\frac{1}{e}\ .
 \end{eqnarray*}
 \qed
 \bigskip

 \noindent
 We now continue with the estimation of~(\ref{Deltasum2}). As a
 result of Lemma~\ref{zeroterms} and Lemma~\ref{eigenvalues}, and
 recalling that $k\!\ge\!n^{\beta'}$,\
 it holds for any $n\ge n_1\mbox{\footnotesize{$\vee$}}\,n_2$\ \,that
 \begin{align} \Label{4.10}
  \left|\dnc\right| 
       & \le \sum_{\la\in\Pn'} e^{-n^{\beta'+\eta-1}}|\lip \gn,
              \Xla\rip|\,|\lip f_n,\Xla\rip| \,.
 \end{align}
 To estimate the number of terms in~(\ref{4.10}), note that the Young diagram
 $\Ul$ of any $\la\in\Pn$ with $B_\la\!=\!s$ consists of an
 $s\mbox{\small{$\times$}}s$ square of cells, with (certainly no more than
 $n\!-\!1$) cells added to each one of the square's $s$ rows and $s$ columns.
 Ignoring the various additional constraints, there are $n^{2s}$
 ways of making such additions, and thus for any $t>0$,\ \
 $\#\{\la\!\in\!\Pn\,:\,B_\la\le t\}\le tn^{2t}$, so that
 $$ \#\Pn'\ \le\ \#\{\la\!\in\!\Pn: B_\la\le n^\beta\}
 \le n^\beta n^{2n^\beta}\le e^{3n^\beta\log n}.$$
 As for the terms in~(\ref{4.10}),\ \ \ $|\lip f_n,\Xla\rip|\le 1$\ \
 by the Cauchy--Schwartz inequality,  and
 \begin{equation*}
  \sup_{\la\in\Pn'}|\lip\gn,\Xla\rip|\le
  \sup_{\underset{\scriptstyle{\gamma\in\Pn:N_\ga\le n^\beta}}
                                  {\la\in \Pn'}}| \chi_\la (\gamma)|
               \le n^{n^\beta}=e^{n^\beta \log n}
 \end{equation*}
 where the second inequality follows from applying \MN's rule at most $n^\beta$
 times, each time with not more that $n$ terms in the sum~(\ref{MN}).

 \noindent
 The above and~\req{4.10} imply that for all $n\!\ge\!\max\{n_0,n_1,n_2\}$
 \[ \left|\dnc\right|\le \exp\left\{-n^\beta
           \left(n^{(\beta'-\beta)-(1-\eta)}-4\log n\right)\right\}\ .\]
 Eventually, thus,\ \  $\left|\dnc\right|\le e^{-\frac{n^\beta}{2}}$,
 which concludes the proof of Theorem~\ref{the-xxx}\ .
 \qed
 \section{Proof of Vershik's conjecture} \Label{Vershik}
  This section is devoted to the proof of Conjecture~\ref{conj-vershik}, which
  we restate as
  \begin{theorem}  \Label{euni}
   If $\mu\in\calM_1(\Omega_1)$ is CCF-invariant then $\mu$ is the
   Poisson--Dirichlet measure $\wmu_1$.
  \end{theorem}
  The main ingredients in its proof have been established in
  Sections~\ref{CCF},~\ref{coupling} and~\ref{DCFconv} and are, respectively,
  the a priori finite moment estimate Proposition~\ref{m1}, the couplings
  with approximating \DCFn's of Theorem~\ref{m7}, and the fast convergence to
  equilibrium of the \DCFn\ chains in the sense of  Theorem~\ref{the-xxx}.
 \bigskip

 \noindent
 {\bf Proof of Theorem~\ref{euni}:}\\
   Let $\Omega_1'=\{p\in\Omega_1\,:\,\exists\ \mbox{infinitely many}\
                  n\in\Nat\ \mbox{such that}\ p_n>\sum_{j>n}p_j\}$.
   We shall show that
   \begin{eqnarray}
    &&\wmu_1(\Omega_1')=1    \Label{support} \\
    &&\{C\,\mbox{\large{$\cap$}}\,\Omega_1'\,:\, C\in\calC\}\subset
                  \calB_{\Omega_1'}\ \ \mbox{is measure determining on\ \ }
                    (\Omega_1',\calB_{\Omega_1'})     \Label{calCisenough}\\
    &&\mu(C)=\wmu_1(C)\hspace{1.2cm} \forall C\in\calC  \Label{okonprime}
   \end{eqnarray}
 which together imply in particular that $\mu(\Omega_1')=1$, and indeed the
 theorem's statement as well.
 \medskip

 \noindent
 {\bf Proof of~(\ref{support}):}\ Recall $\wmu_1$'s description as the law of
 the nonincreasing rearrangement of the uniform stickbreaking process $X_n$
 (with $Y_n=1-\sum_{j<n}X_j$ the remaining stick length prior to the $n$-th
 break and $X_n=U_nY_n$) and define\ $\tau_1\!\!=\!\!1$,\ \
 $\tau_{k+1}\!\! =\min\{n\!>\!\tau_k\ :\ X_n\mbox{\small$\wedge$}
     (Y_n\!-\!X_n)<X_j,\ \forall j\!\le\!\tau_k\,\}$ \ \ for $k\!\!\ge\!\!1$.
 Since a.s. $X_n\mbox{\small$\searrow$}\,0$, each $\tau_k$ is finite.
 We claim that
  \begin{equation} \Label{Ak}
   A_k\!:=\!\{U_{\tau_k}\!>\!\frac{1}{2}\}\!
         =\!\{X_{\tau_k}\!>\!\sum_{j>\tau_k}X_j\} \ \mbox{are independent,}\ \
            \ \ P(A_k)\!=\!\frac{1}{2}\ \ \ \forall k\ .\ \ \ \ \
  \end{equation}
  This implies that a.s.
  $U_{\tau_k}>\frac{1}{2}$ infinitely often, and these $n\!=\!\tau_k$ will
  be the ones alluded to in $\Omega_1'$'s definition. Indeed, on  $A_k$\ \
  \ \  $\sum_{_{j>\tau_k}}\!\!X_j\!<\!X_i\ \ \ \ \forall i\le\tau_k$, so
  that the nondecreasing permutation of the $X_i$'s decouples on $[1,\tau_k]$
  and  $(\tau_k,\infty)$ and thus\ \
  $p_{\tau_k}=\min_{i\le\tau_k}\!X_i>\sum_{_{j>\tau_k}}\!\!X_j
             =\sum_{_{j>\tau_k}}\!p_j.$\\

  \noindent
  To prove~(\ref{Ak}), represent the splitting variables as
  $U_n\!=\!\left\{\begin{array}{ll}
                 V_n   & \mbox{\ if }\eta_n\!=\!1 \\
                 1\!-\!V_n & \mbox{\ if }\eta_n\!=\!0 \\ \end{array}\right.$,\ \
  where $V_n\!\sim\! U[0,1]$ and $\eta_n\!\sim\! \mbox{Bernoulli}(0.5)$ are
  independent of each other, and write
  $A_k\!=(B_k\bigcap C_k)
     \bigcup\,(B_k^{^{\mbox{\tiny C}}}\bigcap C_k^{^{\mbox{\tiny C}}})$,
  with $B_k\!=\!\{V_{_{\tau_k}}\!\!>\!\!\frac{1}{2}\}$
  and $C_k\!=\!\{\eta_{\tau_k}\!\!=\!\!1\}$.
  The $\tau_k$ are $\calF$--stopping times, where
  $\calF_n=\sigma(V_1,\ldots,V_n,\eta_0,\ldots,\eta_{n-1})$
  \ (arbitrarily set $\eta_0\!=\!1$)
  so that $B_k\in\calF_{\tau_k}$ and $C_k$
  is independent of $\calF_{\tau_k}$ (in particular $P(C_k)=0.5$)\ for all $k$.
  For any $D\in\calF_{\tau_k}$
  \begin{eqnarray*}
    P(D\cap A_k)\!\!
      &=&\!\!P((D\cap B_k)\cap C_k)
           +P((D\cap B_k^{^{\mbox{\tiny C}}})\cap C_k^{^{\mbox{\tiny C}}})\\
      &=&\!\!P(D\cap B_k)P(C_k)
           +P(D\cap B_k^{^{\mbox{\tiny C}}})P(C_k^{^{\mbox{\tiny C}}})
        =\frac{1}{2}\,P(D)\ .
  \end{eqnarray*}
  Choosing first\ $D\!=\!\Omega_1$\ and then\ $D\!=\cap_{j\in J}A_j$\ \ with
  $J\subset\{1,\ldots,k-1\}$\ (indeed,
  $A_j\in\calF_{\tau_j+1}\subset\calF_{\tau_k}$ for $j\!<\!k$),\ \
  we  respectively obtain $P(A_k)=0.5$ and the independence of the $A_n$'s.
  We have thus proved~(\ref{Ak}) and thus~(\ref{support}).
  \bigskip

 \noindent
  {\bf Proof of~(\ref{calCisenough}):}\ Fix $\eps\!>\!0,\ p\!\in\!\Omega_1'$,
  and choose $k$ large enough so that
  $0\!<\!q\!:=\!\sum_{_{j>k}}\!p_j\!<\!\min(p_k\,,\,\frac{\eps}{4})$.
  Then let $\delta\!=\!\frac{q\mbox{\scriptsize{$\wedge$}}(p_k-q)}{k+2}$\ \,
  and\ \,$a_i\!=\!p_i\!-\!\delta,\ \ b_i\!=\!p_i\!+\!\delta$\ \ \
  for $i\!=\!1,\ldots,k$. We claim that $(\aaa,\bbb)\!\in\!I_k$. Indeed,
  \begin{eqnarray*}
   \sum_{i\le k}b_i&=&\sum_{i\le k}(p_i+\delta)\le 1-q+\frac{k}{k+2}q<1 ,
                                           \hspace{2cm}\mbox{whereas}\\
   a_k+\sum_{i\le k}a_i&=&(p_k-\delta)+\sum_{i\le k}(p_i-\delta)
                                              =(p_k+1-q)-(k+1)\delta>1\ .
  \end{eqnarray*}
  By definition $p\in\Cab$. Moreover, for any $x\in\Cab$
  \begin{equation}
    |x-p|_1\le 2\sum_{j=1}^k |x_j-p_j|+2\sum_{j=k+1}^\infty p_j
           \le 2k\delta+2q\le 2(1+\frac{k}{k+1})q<\eps\ ,\nonumber
  \end{equation}
  which shows that for any open $l_1$--ball $B_\eps(p)$ in $\Omega_1'$ there is
  some $C\!\in\!\calC$ such that $p\!\in\!C\cap\Omega'_1\subset B_\eps(p)$.
  In other words, $\{C\,\mbox{\large{$\cap$}}\,\Omega_1',\  C\in\calC\}$
  generates $\Omega'_1$'s topology.

  \noindent
  To conclude the proof of~(\ref{calCisenough}) we need to check that $\calC$
  is closed under intersections. For any $j\le k$ then, let
  $(\aaa_1,\bbb_1)\in I_j$  and  $(\aaa_2,\bbb_2)\in I_k$, and if  $j<k$ denote
  ${a_1}_{_i}\!=\!0$\ and\ ${b_1}_{_i}\!=\!1$\ \ for $i\!=\!j+1,\ldots,k$.
  It  follows immediately that $(\aaa,\bbb)$ defined by
  \ $a_i\!=\!{a_1}_{_i}\!\vee{a_2}_{_i}$ \ and
  \ $b_i\!=\!{b_1}_{_i}\!\wedge{b_2}_{_i}$\ \ for\ $i=1,\ldots,k$\ \ belongs
  to $I_k$, and $C_{{\saa_1},{\sbb_1}}\cap C_{{\saa_2},{\sbb_2}}=\Cab$.

 \noindent
 {\bf Proof of~(\ref{okonprime}):}\
First note that if $(\aaa,\bbb)\!\in\!I_k$\ then\
   $((1\!+\!\eps)\aaa,(1\!-\!\eps)\bbb)\!\in\!I_k$\
for all $\eps$ in some neighborhood of $0$,
   and if $C:=\Cab\ne\emptyset$, then so  is
$\Ceps\!:=C_{(1+\eps)\saa,(1-\eps)\sbb}$
   for all $\eps$ in a neighborhood of $0$.

 \noindent
 Once we show that for all $\eps>0$ small enough,
\[  \wmu_1(\Cepsm)\geq \mu(C)\geq \wmu_1(\Ceps)\,,\]
  let $\eps\!\searrow
\!0$ and use  $\wmu_1(\partial C)=0$ to obtain $\mu(C)=\wmu_1(C)$
  for every
$C\!\in\!\calC$,
  thus proving~(\ref{okonprime}) and with it the theorem.

  \noindent
  Let $\frac{2}{5}\!<\!\alpha\!<\!\beta\!<\!\gamma\!<\!\frac{1}{2}$ be three
  otherwise arbitrary numbers. Since by Proposition~\ref{m1}\ $\mu$
  satisfies~(\ref{m7}), we consider for every $n\!\in\!\Nat$ the probability
  measure\ $Q^{(n)}\!=\!Q_\mu^{(n)}$\ introduced in Proposition~\ref{m2}
  which is defined on a  space which supports both a CCF Markov chain
  $p(\cdot)$ with $\mu$ as its stationary marginal and a DCF$\mbox{}^{(n)}$
  Markov chain $\ell^{(n)}(\cdot)$ which ``emulates" $p(\cdot)$ in terms of its
  initial law (cf.~(\ref{m8})) and in the sense that they remain close after
  $n^\gamma$ units of time (cf.~(\ref{m9})). For any $n\in\Nat$,
  \begin{eqnarray}
   \mu(C)-\wmu_1(\Ceps)
      & = &Q^{(n)}(p(\nga)\in C)-\wmu_1(C^{(\eps)}) \nonumber\\
      &\ge&Q^{(n)}(p(\nga)\in C)-Q^{(n)}\left(\mbox{\large{$\frac{1}{n}$}}\,
                              \ell^{(n)}(\nga)\in C^{(\eps)}  \right)\ \nonumber\\
      && \hspace{.5cm}-\ \left|Q^{(n)}\left(\mbox{\large{$\frac{1}{n}$}}\,
             \ell^{(n)}(\nga)\in \Ceps\right) -\pi^{(n)}(n\Ceps)\right|\
          +\left(\pi^{(n)}(n\Ceps)-\wmu_1(\Ceps)\right) \nonumber \\
      & =: &D^{(\eps)}_1-D^{(\eps)}_2+D^{(\eps)}_3.\nonumber
  \end{eqnarray}
  The first term is estimated using a simple union bound with
$\eps':=\eps \min\{a_i\, :\, 1\leq i\leq k\, \}$, and~(\ref{m9}):
  \begin{eqnarray*}
   D^{(\eps)}_1&\ge&-Q^{(n)}\left(\left|p(\nga)-\mbox{\large{$\frac{1}{n}$}}\ell^{(n)}(\nga)
                                                     \right|_1>\eps'\right)
       \ge \frac{-1}{\eps'}E_{Q^{(n)}}\left|p(\nga)
                           -\mbox{\large{$\frac{1}{n}$}}\ell^{(n)}(\nga)\right|_1
            \ \ \underset{n\to\infty}{\mbox{\Large{$\longrightarrow$}}}\ \ 0\ \ .
  \end{eqnarray*}
  To estimate $D^{(\eps)}_2$ we would like to apply Theorem~\ref{the-xxx} to the sequence
  of discrete processes $\ell^{(n)}(\cdot)$. Their initial laws, however, are
  guaranteed by~(\ref{m8}) to be only {\it nearly} supported on
  $\calP_{n,\beta}$, respectively, but not totally as required by
  Theorem~\ref{the-xxx}. Define thus\ \
  $\widetilde{Q}^{(n)}(\cdot)
    :=Q^{(n)}\left(\ \cdot\ |\ \,\ell^{(n)}(0)\!\in\!\calP_{n,\beta}\right)$;\ \ \
  obviously $\widetilde{Q}^{(n)}(\ell^{(n)}(0)\in\calP_{n,\beta})=1$,\ and under
  $\widetilde{Q}^{(n)},\ \ \ell^{(n)}(\cdot)$\ remains a \DCFn\ chain. Then,
  \begin{eqnarray*}
    D^{(\eps)}_2&\le&\left|\widetilde{Q}^{(n)}\left(
             \ell^{(n)}(\nga)\in nC^{(\eps)}\right)-\pi^{(n)}(nC^{(\eps)}) \right|
               +\|\widetilde{Q}^{(n)}-Q^{(n)}\|_{_{_{\mbox{\scriptsize var}}}}
             \ \underset{n\to\infty}{\mbox{\Large{$\longrightarrow$}}}\ \ 0.
  \end{eqnarray*}
   Here we applied Theorem~\ref{the-xxx} for the first term, while
   $ \|\widetilde{Q}^{(n)}-Q^{(n)}\|_{_{_{\mbox{\scriptsize var}}}}\le
     \frac{Q^{(n)}(\ell^{(n)}(0)\not\in\calP_{n,\beta})}
          {Q^{(n)}(\ell^{(n)}(0)\in\calP_{n,\beta})}\to 0$
   by~(\ref{m8}).

   \noindent
   Finally, recall that $\pins(n\cdot)\rightarrow \wmu_1$ weakly
(\cite{Kinew,VS}),
   and since $\Ceps$ satisfies
$\wmu_1(\partial\Ceps)=0$, it follows that
   $\underset{\mbox{\scriptsize{$n\!\!\to\!\!\infty$}}}{{\mbox{$\lim$}}}
   \pins(n\Ceps)= \wmu_1(\Ceps)$. Thus
   $\underset{\mbox{\scriptsize{$n\!\!\to\!\!\infty$}}}{{\mbox{$\lim$}}}
       D^{(\eps)}_3= 0$.\ Consequently
   \[ \mu(C)\!-\!\wmu_1(\Ceps)\ge
       \underset{\mbox{\scriptsize{$n\!\!\to\!\!\infty$}}}{\underline{\mbox{$\lim$}}}D^{(\eps)}_1
      -\underset{\mbox{\scriptsize{$n\!\!\to\!\!\infty$}}}{\overline{\mbox{$\lim$}}}D^{(\eps)}_2
      +\underset{\mbox{\scriptsize{$n\!\!\to\!\!\infty$}}}{\underline{\mbox{$\lim$}}}D^{(\eps)}_3
      \ge 0\ .\]
The reverse inequality $\wmu_1(\Cepsm)\geq \mu(C)$ is obtained similarly from
\[\wmu_1(\Cepsm)-\mu(C)\geq -D_1^{(-\eps)}-D_2^{(-\eps)}-D_3^{(-\eps)}.\]
\qed
 \section{Appendix}
  \paragraph{Proof of Proposition~\ref{m1}}
Consider the partition of $(0,1]$ by
$J_n:=(2^{-n-1},2^{-n}]\quad(n\geq 0)$ and define on $\Omega_1$ the random
variables
\[W_n:=\sum_{i\geq 1}p_i\won_{p_i>2^{-n}}\quad(n\geq 1).\]
Fix $n\geq 1$. If two intervals are merged then $W_n$ can only
increase and if some interval is split then  $W_n$ can only
decrease. We call the increment in the case of merging
$\Delta_+\geq 0$ and the loss in the case of splitting
$\Delta_-\geq 0$.
 Given $p$, we can bound $\Delta_+$  by
\begin{eqnarray*}
\Delta_+&\geq& \sum_{i\ne j}p_i^2p_j\won_{p_i\in
J_n}\won_{p_j>2^{-n-1}}\\ &=& \left(\sum_{i}p_i^2\won_{p_i\in
J_n}\right) \left(\sum_j p_j\won_{p_j>2^{-n-1}}\right)-
\sum_{i}p_i^3\won_{p_i\in J_n}\\ &\geq&
2^{-2n-2}\left(\sum_{i}\won_{p_i\in J_n}\right)\left(\sum_j
p_j\won_{p_j>2^{-n-1}}\right) -2^{-2n}\sum_{i}p_i\won_{p_i\in J_n},
\end{eqnarray*}
and compute $\Delta_-$ as
\begin{eqnarray*}
\Delta_-&=& \sum_ip_i^2\int_0^1xp_i\won_{xp_i\leq 2^{-n}<p_i}+
(1-x)p_i\won_{(1-x)p_i\leq 2^{-n}<p_i}\ dx\\ &=&2 \sum_ip_i^3
\won_{2^{-n}<p_i}\int_0^1 x \won_{x\leq 2^{-n}/p_i}\ dx\ =\ 2
\sum_ip_i^3
\won_{2^{-n}<p_i}\left[\frac{x^2}{2}\right]_0^{2^{-n}/p_i}\ =\
2^{-2n}\sum_i p_i \won_{2^{-n}<p_i}.
\end{eqnarray*}
Therefore,
\begin{eqnarray*}
\Delta_+-\Delta_-&\geq& 2^{-2n}\left(
\frac{1}{4}\left(\sum_{i}\won_{p_i\in J_n}\right)\left(\sum_j
p_j\won_{p_j>2^{-n-1}}\right)- \sum_i p_i \won_{2^{-n-1}<p_i}\right)\\ &
\geq&
 2^{-2n}\left(\frac{1}{4}
\left(\sum_{i}\won_{p_i\in J_n}\right)\left(\sum_j
p_j\won_{p_j>2^{-n-1}}\right)- 1\right).
\end{eqnarray*}
Since  $\int\Delta_+-\Delta_-\ d\mu=0$ due to stationarity this
implies
\begin{eqnarray}
\Label{math51} 4&\geq&\int \left(\sum_{i}\won_{p_i\in
J_n}\right)\left(\sum_j p_j\won_{p_j>2^{-n-1}}\right)\ d\mu\\ &\geq
&2^{-n-1}\int\left(\sum_{i}\won_{p_i\in J_n}\right)^2\ d\mu\ \geq\
2^{-n-1}\left(\int\sum_{i}\won_{p_i\in J_n}\ d\mu\right)^2.\nonumber
\end{eqnarray}\abel{math51}
Since this holds for any $n\geq 1$ we get
\begin{equation}\Label{zovim}
\int \sum_{i}\won_{p_i\in J_n}\ d\mu= O(2^{\beta
n})\mbox{\quad($n\to\infty$)}
\end{equation}\abel{zovim}
with $\beta=1/2$. Therefore, for any $\alpha>\beta=1/2$,
\begin{equation}\Label{har}
\int\sum_{i}p_i^\alpha\ d\mu\leq \sum_{n\geq 0}2^{-\alpha n}
\int\sum_{i} \won_{p_i\in J_n}\ d\mu \leq c  \sum_{n\geq
0}2^{n(\beta-\al)}<\infty,
\end{equation}\abel{har}
for some constant $c>0$, thus proving (\ref{hi}) for
$\al>1/2$. We shall now use this result to extend it to all
$2/5<\al<1$, as required. To this end, observe that we have due to
(\ref{math51}) for arbitrary $0<\beta<1$,
\begin{eqnarray}
4&\geq & \int \left(\sum_{i}\won_{p_i\in
J_n}\right)\left(\sum_j p_j\won_{p_j>2^{-n-1}}\right) \won\left\{\sum_j
p_j\won_{p_j>2^{-n-1}}>2^{-n\beta}\right\}\ d\mu\nonumber\\ &\geq&
2^{-n\beta} \int\sum_{i}\won_{p_i\in J_n} \won\left\{\sum_j
p_j\won_{p_j>2^{-n-1}}>2^{-n\beta}\right\}\ d\mu\nonumber\\ &=&
2^{-n\beta}\left( \int\sum_{i}\won_{p_i\in J_n}\ d\mu-\int
\sum_{i}\won_{p_i\in J_n} \won\left\{\sum_j p_j\won_{p_j>2^{-n-1}}\leq
2^{-n\beta}\right\}\ d\mu\right)\nonumber\\ &\geq&
2^{-n\beta}\left( \int\sum_{i}\won_{p_i\in J_n}\
d\mu-\frac{2^{-n\beta}}{2^{-n-1}}\mu\left[ \sum_j
p_j\won_{p_j>2^{-n-1}}\leq 2^{-n\beta}\right]\right)\nonumber\\
&\geq& 2^{-n\beta}\left( \int\sum_{i}\won_{p_i\in J_n}\
d\mu-2^{-n\beta +n+1}\mu[\forall i: p_i\leq 2^{-n\beta}]\right).
\Label{fs}
\end{eqnarray}\abel{fs}
To bound the $\mu$-probability in the last expression we recall
from
 (\ref{har}) that for $1/2<\al<1$,
\begin{equation}\Label{mich}
\infty> \int\sum_{i}p_i^\alpha\ d\mu\geq
\int\left(\sum_{i}p_i^\alpha\right) \won\left\{\forall i: p_i\leq
2^{-n\beta}\right\}\ d\mu
\end{equation}\abel{mich}
for any $n\geq 0$. On the event $\{\forall i: p_i\leq
2^{-n\beta}\}$, by Jensen's inequality,
\[\sum_i p_i^\al=\sum_i p_i\cdot p_i^{\al-1}\geq
\left(\sum_i p_i\cdot p_i\right)^{\al-1}\geq 2^{n\beta(1-\al)}.\]
Therefore, for all $1/2<\al<1$ due to ({\ref{mich}}), $\mu[\forall
i: p_i\leq 2^{-n\beta}]=O(2^{-n\beta(1-\al)})$ as $n\to\infty$.
 Substituting this into (\ref{fs}) we get that for any $1/2<\al<1$
\[\int \sum_{i}\won_{p_i\in J_n}\ d\mu= O(2^{n\beta}+ 2^{(-n\beta +n)
-n\beta(1-\al)})= O(2^{n\max\{\beta, 1-\beta(2-\al)\}}).\] The
choice $\beta=(3-\al)^{-1}$ minimizes $\max\{\beta,
1-\beta(2-\al)\}$ and therefore yields that (\ref{zovim}) and
consequently also (\ref{har}) and (\ref{hi}) hold for any
$\alpha,\beta>(3-1/2)^{-1}=2/5$. \qed

 \noindent
 \remark A posteriori, once it has been established that $\mu$ must be
 the Poisson-Dirichlet law, Proposition~\ref{m1} holds for
 all $\alpha>0$ since by~\cite[(20)]{Ki1}
 \[ \int\sum_{i\ge 1}p_i^\alpha d\wmu_1=\int_0^1 x^{\alpha-1}\,dx
         =\frac{1}{\alpha}\ .\]

\noindent
{\bf Note added in proof:}
In [13], the question was raised as to whether the state 
$s:=(1,0,\ldots)\in \Omega_1$ is recurrent for the CCF. Our techniques allow
one to respond affirmatively to this question. Indeed, let $X^{(n)}(k)$ denote
the state, at time $k$, of a DCF${}^{(n)}$ initialized at $X^{(n)}(0)=(n,0,\ldots)=:s_n$.
The recurrence of $s$ for the CCF then follows, by the coupling
introduced in Theorem 3.1,  from the existence of a constant $C$ independent of $n$ and
$k<n$ such that
$P^{(n)}
(X^{(n)}(2k)=s_n)\geq C/k$. To see the last estimate, note that by the character decomposition
at the beginning of the proof of Theorem 4.1, it holds that
$$
P^{(n)}(X^{(n)}(2k)=s_n)
     =\sum_{\lambda\in
{{\cal P}_n}}{\theta_\lambda^{(n)^{2k}}}
                \langle (\pi_S^{(n)}(s_n))^{-1}{\bf 1}_{s_n},{\chi_{_\lambda}}\rangle
\langle{\bf 1}_{s_n} ,{\chi_{_\lambda}}\rangle
=
     \sum_{\lambda\in
{{\cal P}_n}}{\theta_\lambda^{(n)^{2k}}}
\pi_S^{(n)}(s_n) {\chi_{_\lambda}}^2(s_n)
\,.$$
From (1.1), $\pi_S^{(n)}(s_n)=1/n$.
By the Murnaghan--Nakayama rule, $\chi_{_\lambda}(s_n)=0$ unless $\lambda=(i,1,1,\ldots)$
for some $i\in \{1,\ldots,n\}$, in which case $|\chi_{_\lambda}(s_n)|=1$.  Using (4.6),
one has that for such $\lambda$, $\theta_\lambda^{(n)}=(2i-n-1)/(n-1)$. 
Thus,
$$
P^{(n)}(X^{(n)}(2k)=s_n)=\frac1n \sum_{i=1}^n \left(1-\frac{2i}{n-1}\right)^{2k}
\geq 
\frac1n \sum_{i=1}^{n/4} \left(1-\frac{2i}{n-1}\right)^{2k}
\geq
\frac1{2n} \sum_{i=1}^{n/4} e^{-4ki/(n-1)}\,,$$
where we used that $(1-x)\geq e^{-x}/2$ for $x<1/\sqrt{2}$. This yields the claim.
%

\vspace{1cm}

\begin{tabular}{ll}
Persi Diaconis &Eddy Mayer-Wolf \\
   Dept.\ of
Mathematics &  Dept.\ of Mathematics\\
and Dept.\  of Statistics& Technion\\
Stanford University &
  Haifa 32000, Israel\\
Stanford, CA 94305, U.S.A.
& emw@tx.technion.ac.il\\
 &\\
Ofer Zeitouni& Martin Zerner\\
  Departments of Math. and of EE& Dept.\ of Mathematics\\
Technion, Haifa 32000, Israel&
Stanford University  \\
and Dept.\ of Math.&
Stanford, CA 94305, U.S.A. \\
University of Minnesota&\\
Minneapolis, MN 55455, USA&\\
zeitouni@math.umn.edu&zerner@math.stanford.edu\\
www-ee.technion.ac.il/$\!{}_{\mbox{\~{}}}$zeitouni
&
math.stanford.edu/$\!{}_{\mbox{\~{}}}$zerner
\end{tabular}


\begin{thebibliography}{99}
 \bibitem{arratia} R.\ Arratia, A.\ D.\ Barbour and S.\ Tavar\'{e},
   {\it Logarithmic Combinatorial Structures: a Probabilistic Approach},book,
   preprint (2002), http://www-hto.usc.edu/books/tavare/ABT/index.html .
 \bibitem{brooksmakover1} R.\ Brooks  and E.\ Makover, Random construction of
   Riemann surfaces, \textit{Preprint}, Department of Mathematics, Technion
   (1997).
 \bibitem{brooksmakover0} R.\ Brooks and E.\ Makover, Belyi surfaces, in:
   \textit{Entire functions in modern analysis. Israel Math. Conf. Proc.},
   {\bf 15}, pp.~37-46, Bar Ilan Univ. (2001)
 \bibitem{diaconis} P.\ Diaconis, {\it Group Representations in Probability
   and Statistics}, Institute of Mathematical Statistics, Hayward,
   Lecture Notes--Monograph Series, {\bf 11} (1988).
 \bibitem{diaconishanlon} P.\ Diaconis and P.\ Hanlon, Eigen analysis for
   some examples of the Metropolis algorithm, in: {\it
   Hypergeometric functions on domains of positivity, Jack polynomials, and
   applications. Contemp.\ Math.} {\bf 138},  pp.~99-117, AMS, Providence (1992).
 \bibitem{flattoetal} L.\ Flatto, A.M.\ Odlyzko and D.B.\ Wales, Random Shuffles
   and Group Representations, {\it Ann. Prob.} {\bf 13}, (1985) pp.~154--178.
 \bibitem{gamburdmakover} A.\ Gamburd  and E.\ Makover, On the genus of a
   random Riemann surface, \textit{Contemp.\ Math.} {\bf 311} (2002), pp.\
   133-140.
 \bibitem{gnedin} A.\ Gnedin and S.\ Kerov, A characterization of GEM
   distributions, {\it Combin. Probab. Comp.} {\bf 10},  (2001),
   pp.~213--217.
 \bibitem{gnedin1} A.\ Gnedin and S.\ Kerov, Fibonacci solitaire,
{\it Random Struct. Alg.} {\bf 20} (2002), pp.~71-88.
 \bibitem{jameskerber} G.\ James and A.\ Kerber, {\it The Representation Theory
   of the Symmetric Group}, Addison--Wesley, Reading, Massachusetts (1981).
 \bibitem{Ki1} J.\ F.\ C.\ Kingman, Random discrete distributions,
   {\it J. Roy. Statist. Soc. Ser. B} {\bf 37} (1975), pp. 1--22.
 \bibitem{Ki2} J.\ F.\ C.\  Kingman, {\it Poisson Processes}, Oxford University
   Press, Oxford (1993).
\bibitem{Kinew} J.\ F.\ C.\ Kingman, The population structure associated with
the Ewens sampling formula, {\it Theoretical Population Biology} {\bf 11}
(1977), pp. 274--283.
 \bibitem{MZZ}E.\ Mayer-Wolf, O.\ Zeitouni and M.P.W.\ Zerner, Asymptotics of
   certain coagulation-fragmentation processes and invariant Poisson-Dirichlet
   measures, {\it Electr.\ J.\ Prob.} {\bf 7} (2002), paper no.\ 8, pp.~1--25
 \bibitem{pitman} J.~Pitman, Poisson--Dirichlet and GEM invariant distributions
   for split-and-merge transformations of an interval partition, {\it Combin.
   Prob. Comp.}, to appear.
 \bibitem{Tsilevich} N.\ V.\ Tsilevich, Stationary random partitions of
   positive integers, {\it Theor. Probab. Appl.} {\bf 44} (2000), pp. 60--74.
 \bibitem{Tsilevich1} N.\ V.\ Tsilevich, On the simplest split-merge operator
   on the infinite-dimensional simplex, {\it PDMI PREPRINT 03/2001}, (2001).
   ftp://ftp.pdmi.ras.ru/pub/publicat/preprint/2001/03-01.ps.gz
 \bibitem{VS} A. M. Vershik and A. A. Shmidt, Limit theorems arising in the
   asymptotic theory of symmetric groups, I., {\it Th. Prob. Appl.} {\bf 22}
   (1977), pp. 70--85.
 \end{thebibliography}
\end{document}